\documentclass[a4paper]{amsart}
\usepackage{amsmath,amssymb,amsfonts}
\usepackage{colordvi}
\usepackage[all]{xy}

\def\mpoint{\;.}
\def\mvirg{\;,}

\def\sou{\underline}

\def\mpn{\medskip\par\noindent}

\def\mmpn{\vskip 1em minus 1em\par\noindent}

\def\sp{\bigskip\par}

\def\smp{\smallskip\par}

\def\CB{{\mathcal B}}
\def\CC{{\mathcal C}}

\def\CE{{\mathcal E}}
\def\CF{{\mathcal F}}

\def\CP{{\mathcal P}}

\def\CR{{\mathcal R}}

\def\Ker{\operatorname{Ker}\nolimits}

\def\Id{\operatorname{id}\nolimits}

\def\Mod{\operatorname{Mod}\nolimits}

\def\Aut{\operatorname{Aut}\nolimits}

\def\Irr{\operatorname{Irr}\nolimits}

\def\tot{\operatorname{tot}\nolimits}
\def\equi{{\scriptstyle{\preceq\succeq}}}

\newcommand{\vdashER}[1]{\mathop{\vdash}_{\scriptscriptstyle #1}\limits}
\newcommand{\nvdashER}[1]{\mathop{\nvdash}_{\scriptscriptstyle #1}\limits}

\def\op{^{op}}
\def\dual{^{\scriptscriptstyle\natural}}

\def\ls#1#2{{\,^{#1}\!#2}}
\def\Iup{I^{\uparrow}}
\def\Idown{I_{\downarrow}}
\def\meet{\wedge}
\def\pointplein{\makebox[0ex]{$\bullet$}}
\def\pointcreux{\makebox[0ex]{$\circ$}}

\def\Z{\mathbb{Z}}
\def\N{\mathbb{N}}

\def\S{\mathbb{S}}

\newcommand{\sumb}[2]{\sum_{{\scriptstyle #1}\atop {\scriptstyle #2}}}

\newcommand{\edge}[2]{\xymatrix{#1\ar@{->-}[r]&#2}}
\def\marc[#1]{\ar@{-}[#1]|(.4){\object@{<}}}
\def\mard[#1]{\ar@{-}[#1]|(.5){\object@{>}}}
\def\marb[#1]{\ar@{-}[#1]|{\object+{  }}}
\newcommand{\fleche}[2]{\xymatrix@C=4ex{*!U(0.2){#1\;}&*!U(0.5){\;#2}\marc[l]}}
\newcommand{\flecheb}[2]{\xymatrix@C=4ex{*!U(0.2){#1\;}&*!U(0.1){\;#2}\marc[l]}}

\def\pf{\par\bigskip\noindent{\bf Proof~: }}
\def\endpf{~\hfill\rlap{\hspace{-1ex}\raisebox{.5ex}{\framebox[1ex]{}}\sp}\bigskip\pagebreak[3]}

\makeatletter
\renewenvironment{enumerate}{\ifnum \@enumdepth >3 \@toodeep\else
       \advance\@enumdepth \@ne
       \edef\@enumctr{enum\romannumeral\the\@enumdepth}\list
       {\csname  label\@enumctr\endcsname}{\setlength{\topsep}{1ex}
 \setlength{\itemsep}{0 pt}\usecounter
         {\@enumctr}\def\makelabel##1{\hss\llap{##1}}}\fi}{\endlist}

\def\@seccntformat#1{\csname the#1\endcsname.\quad}

\def\section{\pagebreak[3]\setcounter{prop}{0}\setcounter{equation}{0}\@startsection{section}{1}{\z@}{4ex plus  6ex}{2ex}{\center\reset@font \large\bf}}
\newcommand{\subsect}[1]{\medskip\par\noindent\pagebreak[3]\refstepcounter{subsection}\refstepcounter{prop}{\bf \thesection.\arabic{prop}.\ #1.\ }}
\makeatother

\def\theprop{\thesection.\arabic{prop}}

\renewenvironment{equation}{\refstepcounter{subsection}\refstepcounter
{prop}$$}{\leqno{\bf (\theprop)}$$}

\newenvironment{enonce}[1]{\pagebreak[3]\refstepcounter{prop}\mmpn
{{\bf  \thesection.\arabic{prop}.\ #1.}}\begin{it} }{\end{it}\smp}
\def\thesection{\arabic{section}}

\newcommand{\result}[1]{\begin{enonce}{#1}}
\newcommand{\fresult}{\end{enonce}}

\newcommand{\mbigvee}[1]{\mathop{\bigvee}_{#1}\limits}

\newcommand{\mbigwedge}[1]{\mathop{\bigwedge}_{#1}\limits}

\begin{document}

\title[Boolean matrices, correspondence functors, and simplicity]
{The algebra of Boolean matrices, correspondence functors, and simplicity} 

\author{Serge Bouc}
\author{Jacques Th\'evenaz}
\date\today



\begin{abstract}
We determine the dimension of every simple module for the algebra of the monoid of all relations on a finite set (i.e. Boolean matrices).
This is in fact the same question as the determination of the dimension of every evaluation of a simple correspondence functor.
The method uses the theory of such functors developed in~\cite{BT2,BT3}, as well as some new ingredients in the theory of finite lattices.
\end{abstract}

\maketitle


\section{Introduction}

\medskip
\noindent
Let $k$ be a field and let $\CR_X$ be the $k$-algebra of the monoid of all relations on a finite set~$X$ (also known as Boolean matrices).
This is an algebra of dimension $2^{n^2}$, where $n=|X|$, hence growing very fast in terms of~$n$.
It was considered many years ago in~\cite{CP, Ki, KR, PW, Sc1, Sc2} and more recently in~\cite{BE, Br, Di}, but the dimensions of the irreducible representations of~$\CR_X$ remained unknown in general. Some related work \cite{St1, St2} shows that this domain of research is still active.\par

We solve here the open problem of describing all simple $\CR_X$-modules and finding their dimension. This requires to embed the category of $\CR_X$-modules into the larger category of correspondence functors, namely functors from the category of finite sets and correspondences to the category $k\text{-\!}\Mod$. We use methods of the representation theory of categories, as well as some new ingredients in the theory of finite lattices. The proof is based on very delicate arguments about a system of linear equations which was introduced in~\cite{BT3}. We also deduce a formula for the dimension of the Jacobson radical of~$\CR_X$ (in characteristic zero). The formulas behave exponentially with respect to~$n$.\par

In a previous work~\cite{BT1}, we described all simple modules for the algebra $\CE_X$ of essential relations on~$X$, which is a quotient of~$\CR_X$, but it was then not clear how to extend this result. It is not too hard to show (and known to some specialists) that the simple modules for~$\CR_X$ are classified by isomorphism classes of triples~$(E,R,V)$, where $E$ is finite set with $|E|\leq|X|$, $R$ is a partial order relation on~$E$, and $V$ is a simple module for the group algebra $k\Aut(E,R)$. When $E=X$, we recover the simple modules for the essential algebra~$\CE_X$, but the more difficult cases occur when $|E|<|X|$.\par

Apart from~\cite{BT1}, the main ingredients for this work are our papers~\cite{BT2,BT3} about correspondence functors. It is known that the evaluation at a finite set~$X$ of a simple functor is either zero or a simple module for the algebra~$\CR_X$. Conversely every simple module for the algebra~$\CR_X$ occurs as the evaluation at~$X$ of a simple functor. This provides a way to handle simple modules for the algebra~$\CR_X$ by studying simple correspondence functors. It is this embedding in the larger category of correspondence functors which allows us to prove our results.
A first step, which is not very hard and explained in~\cite{BT2}, is the description of the parametrization of simple correspondence functors $S_{E,R,V}$ by isomorphism classes of triples $(E,R,V)$, where $(E,R)$ is a finite poset (i.e. $R$ is a partial order relation on a finite set~$E$) and $V$ is a simple module for the group algebra $k\Aut(E,R)$.\par

Some fundamental modules and functors play a crucial role in our approach.
Here $k$ is allowed to be an arbitrary commutative ring.
For any finite poset $(E,R)$, we described in~\cite{BT1} a fundamental module $\CP_E f_R$ for~$\CR_E$,
where $\CP_E$ is a quotient algebra of the essential algebra~$\CE_E$, hence also a quotient of~$\CR_E$, while $f_R$ is a suitable idempotent of~$\CP_E$ depending on the order relation~$R$.
From this, we constructed and studied in~\cite{BT2,BT3} a fundamental functor $\S_{E,R}$, which is the key for understanding simple correspondence functors because the simple functor $S_{E,R,V}$ appears as a suitable quotient of the fundamental functor~$\S_{E,R}$.\par

Another main ingredient is the link between correspondence functors and the theory of finite lattices, see~\cite{BT3}.
Associated to any finite lattice~$T$, there is a correspondence functor $F_T$ and a surjective morphism
$$\Theta: F_T\longrightarrow \S_{E,R\op}$$
where $(E,R)$ denotes the full subposet of join-irreducible elements of~$T$ and $R\op$ denotes the opposite relation.
The main problem is to describe the kernel of~$\Theta$ and this gives rise to a complicated system of linear equations which was introduced in~\cite{BT3}.
One of the main contributions of the present paper is to solve this system.
From this solution, a $k$-basis can be found for each evaluation $\S_{E,R\op}(X)$ of a fundamental functor.
Generators are found in Section~\ref{Section-generators} and they are proved to be $k$-linearly independent in Section~\ref{Section-independence}.\par

Turning to simple functors (assuming again that $k$ is a field), we need to pass to a quotient of~$\S_{E,R}$ in order to obtain the simple functor $S_{E,R,V}$.
This requires to show that each evaluation $\S_{E,R}(X)$ has a free right $k\Aut(E,R)$-module structure and that the simple functor $S_{E,R,V}$ is isomorphic to a tensor product $\S_{E,R}\otimes_{k\Aut(E,R)}V$. This nontrivial part of the argument requires the whole of Section~\ref{Section-simple-tensor} and culminates with an explicit formula for the dimension of each evaluation $S_{E,R,V}(X)$ of a simple correspondence functor. The final step, which is easy and explained in Section~\ref{Section-algebra}, is to go back to the algebra~$\CR_X$
and deduce the dimension of every simple $\CR_X$-module, as well as a description of the action of relations on it.\par

There is a classical approach of the classification of simple modules for the algebra of a finite semigroup, going back to the work of Munn and Ponizovsky, using Green's theory of $J$-classes (see the textbook~\cite{CP}, or the more recent article \cite{GMS} for a modern point of view). For the algebra $\CR_X$ we are interested in, we do not use this point of view here for two reasons. First, our approach of the parametrization of simple modules for~$\CR_X$ is not classical, for it is based in an important way on the fundamental module~$\CP_E f_R$ associated to a poset $(E,R)$. Secondly, taking advantage of the link with the theory of correspondence functors and using the functor~$F_T$ associated to a finite lattice~$T$, our main task is the study of the above morphism~$\Theta$ (which is itself based on the fundamental module~$\CP_E f_R$). This problem has no obvious connection with the classical approach to the algebra of a semigroup.
We leave to the interested reader the open question of investigating if it is possible to translate our results in the language of $J$-classes and other concepts used in the classical approach.


\section{Preliminaries on lattices} \label{Section-lattices}

\medskip
\noindent
In this section, we define, in any finite lattice, two operations $r^\infty$ and $\sigma^\infty$, as well as a subset $\widehat G$ of special elements, each lying at the bottom of a totally ordered subset with strong properties. We then prove some results which will play a crucial role in the description of the evaluation of fundamental functors and simple functors.\par

Let us first fix some notation.
By an {\em order} $R$ on a finite set~$E$, we mean a partial order relation on~$E$.
In other words, $(E,R)$ is a finite poset.
We write $\leq_R$ for the order relation, so that $(a,b)\in R$ if and only if $a\leq_R b$.
Moreover $a<_R b$ means that $a\leq_R b$ and $a \neq b$.
The opposite relation $R\op$ of~$R$ is defined by the property that $(a,b)\in R\op$ if and only if $(b,a)\in R$.\par

If $T$ is a finite lattice, we write $\leq_T$, or sometimes simply $\leq$, for the order relation, $\vee$ for the join (least upper bound), $\wedge$ for the meet (greatest lower bound), $\hat{0}$ for the least element and $\hat{1}$ for the greatest element.

\result{Notation and definitions} 
\begin{enumerate}
\item If $(E,R)$ is a finite poset and $a,b\in E$ with $a\leq_R b$, we define intervals
$$\begin{array}{ll}
\qquad[a,b]_E:=\{x\in E\mid a\leq_R x\leq_R b\} \,,\qquad & {]a,b[_E}:=\{x\in E\mid a<_R x<_R b\} \,, \\
\qquad{[a,b[_E}:=\{x\in E\mid a\leq_R x<_R b\} \,,\qquad & {]a,b]_E}:=\{x\in E\mid a<_R x\leq_R b\} \,, \\
\qquad{[a,\cdot[_E}:=\{x\in E\mid a\leq_R x\} \,,\qquad & {]\cdot,b]_E}:=\{x\in E\mid x\leq_R b\} \mpoint
\end{array}$$
When the context is clear, we write ${[a,b]}$ instead of~${[a,b]_E}$.
\item If $T$ is a finite lattice an element $e\in T$ is called {\em join-irreducible}, or simply {\em irreducible},
if, whenever $e=\mbigvee{a\in A} a$ for some subset $A$ of~$T$, then $e\in A$.
In case $A=\emptyset$, the join is~$\hat0$ and it follows that $\hat0$ is not irreducible.
If $e\neq\hat0$ is irreducible and $e=s\vee t$ with $s,t\in T$, then either $e=s$ or $e=t$.
In other words, if $e\neq\hat0$, then $e$ is irreducible if and only if $[\hat0,e[$ has a unique maximal element.
\item If $(E,R)$ is a subposet of a finite lattice~$T$, we say that $(E,R)$ is a {\em full\/} subposet of~$T$ if for all $e,f\in E$ we have~:
$$e\leq_R f \Longleftrightarrow e\leq_T f \mpoint$$
In particular $\Irr(T)$ denotes the full subposet of irreducible elements of~$T$.
\item If $(E,R)$ is a finite poset, $\Idown(E,R)$ denotes the set of {\em lower $R$-ideals} of~$E$, that is, the subsets $A$ of~$E$ such that, whenever $a\in A$ and $x\leq a$, then $x\in A$.
Clearly $\Idown(E,R)$, ordered by inclusion of subsets, is a lattice,
the join operation being union of subsets, and the meet operation being intersection.
Similarly, $\Iup(E,R)$ denotes the set of upper $R$-ideals of $E$, which is also a lattice.
Obviously $\Iup(E,R)=\Idown(E,R\op)$.
\end{enumerate}
\fresult

Note that if $(E,R)$ is the poset of irreducible elements in a finite lattice $T$, then $T$ is {\em generated\/} by~$E$ in the sense that any element $x\in T$ is a join of elements of~$E$.
To see this, define the height of $t\in T$ to be the maximal length of a chain in $[\hat0,t]_T$.
If $x$ is not irreducible and $x\neq\hat0$, then $x=t_1\vee t_2$ with $t_1$ and $t_2$ of smaller height than~$x$.
By induction on the height, both $t_1$ and $t_2$ are joins of elements of~$E$.
Therefore $x=t_1\vee t_2$ is also a join of elements of~$E$.\par

\result{Lemma} \label{principal} Let $(E,R)$ be a finite poset.
\begin{enumerate}
\item The irreducible elements in the lattice~$\Idown(E,R)$ are the lower ideals $]\cdot,e]_E$, where $e\in E$.
Thus the poset $(E,R)$ is isomorphic to the poset of all irreducible elements in~$\Idown(E,R)$
by mapping $e\in E$ to the ideal $]\cdot,e]_E$.
\item $\Idown(E,R)$ is a distributive lattice.
\item For any finite lattice $T$ with $\Irr(T)=(E,R)$, there is a join-preserving surjective map
$f: \Idown(E,R) \longrightarrow T$ which sends any lower ideal $A\in \Idown(E,R)$
to the join $\bigvee_{e\in A}e$ in~$T$.
\item The map $f: \Idown(E,R) \longrightarrow T$ above is bijective if and only if $T$ is a distributive lattice. In that case, $f$ is an isomorphism of lattices.
\end{enumerate}
\fresult

\pf This is not difficult and well-known.
For details, see Theorem~3.4.1 and Proposition~3.4.2 in~\cite{Sta}, and also Theorem~6.2 in~\cite{Ro}.
\endpf

Whenever we use the lattice $\Idown(E,R)$, we shall (abusively) identify $E$ with its image
via the map
$$E\longrightarrow\Idown(E,R) \,,\qquad  e \mapsto \;]\cdot,e]_E \mpoint$$
Thus we view $(E,R)$ as a full subposet of~$\Idown(E,R)$.\par

\result{Notation} \label{r(t)}
Let $T$ be a finite lattice and $(E,R)=\Irr(T)$.
If $t\in T$, then $r(t)$ denotes the join of all elements strictly smaller than~$t$~:
$$r(t)=\bigvee_{a\in [\hat0,t[} a \mpoint$$
It follows that $r(t)=t$ if and only if $t\notin E$.
More precisely, if $t\notin E$ and $t\neq\hat 0$, then $t$ can be written as the join of two smaller elements, so $r(t)=t$,
while if $e\in E$, then $r(e)$ is the unique maximal element of~$[\hat0,e[$.
We put $r^k(t)=r(r^{k-1}(t))$ and $r^\infty(t)=r^n(t)$ if $n$ is such that $r^n(t)=r^{n+1}(t)$.
\fresult

\result{Lemma} \label{rinfty} Let $T$ be a finite lattice, let $(E,R)=\Irr(T)$, and let $t\in T$.
\begin{enumerate}
\item $r^\infty(t)\notin E$.
\item $r^\infty(t)=t$ if and only if $t\in T-E$.
\item If $e\in E$, $r^\infty(e)$ is the unique greatest element of $T-E$ smaller than~$e$.
\item If $t'\in T$ with $t\leq t'$, then $r^\infty(t)\leq r^\infty(t')$.
\item The map $r:T\to T$ is order-preserving.
\end{enumerate}
\fresult

\pf
The proof is a straightforward consequence of the definitions.
\endpf

\result{Lemma} \label{interval} Let $T$ be a finite lattice, let $(E,R)=\Irr(T)$, and let $e\in E$.
Let $n$ be the smallest integer such that $r^n(e)=r^\infty(e)$.
\begin{enumerate}
\item $[r^\infty(e),e]$ is totally ordered and $[r^\infty(e),e]= \{r^n(e), \ldots, r^1(e),e \}$.
\item $]r^\infty(e),e]$ is contained in~$E$.
\item $r^\infty(r^i(e))=r^\infty(e)$ for all $0\leq i\leq n-1$.
\item $[\hat0,e] = \; [\hat0,r^\infty(e)] \;\sqcup\; ]r^\infty(e),e]$.
\end{enumerate}
\fresult

\pf
Since $e\in E$, $r(e)$ is the unique maximal element of~$[\hat0,e[$.
Inductively, $r^i(e)\in E$ for each $0\leq i\leq n-1$ and $r^{i+1}(e)$ is the unique maximal element of~$[\hat0,r^i(e)[$.
It follows that $[r^\infty(e),e]$ is totally ordered and consists of the elements $r^n(e), \ldots, r^1(e), e$.
This proves~(a), (b) and (c).\par

Now let $f\in [\hat0,e]$. Then $f\vee r^\infty(e)\in[r^\infty(e),e]$.
If $f\vee r^\infty(e)=r^\infty(e)$, then $f\in \;[\hat0,r^\infty(e)]$.
Otherwise, $f\vee r^\infty(e)\in \; ]r^\infty(e),e]$, hence $f\vee r^\infty(e)\in E$ by~(b),
that is, $f\vee r^\infty(e)$ is irreducible. It follows that $f\vee r^\infty(e)=f$ or $f\vee r^\infty(e)=r^\infty(e)$.
But the second case is impossible because $f\vee r^\infty(e)> r^\infty(e)$.
Therefore $f\vee r^\infty(e)=f$, that is, $f\in \; ]r^\infty(e),e]$.
\endpf

\result{Notation and definitions} \label{notation-G} Let $T$ be a finite lattice and $(E,R)=\Irr(T)$.
\begin{enumerate}
\item Define $\Lambda E$ to be the subset of~$T$ consisting of all meets of elements of~$E$,
that is, elements of the form $\bigwedge\limits_{i\in I}e_i$ where $I$ is a finite set of indices and $e_i\in E$ for every $i\in I$.
Note that we include the possibility that $I$ be the empty set,
in which case one gets the unique greatest element~$\hat1$.
\item If $t\in T$, define $\sigma(t)$ to be the meet of all the irreducible elements of~$T$ which are strictly larger than~$t$.
Inductively, $\sigma^k(t)=\sigma(\sigma^{k-1}(t))$ and $\sigma^\infty(t)=\sigma^n(t)$ where $n$ is such that $\sigma^n(t)=\sigma^{n+1}(t)$.
\item Define the subset $\widehat G=\widehat G(T)$ to be the set of all elements $t\in T$ such that
$t\notin \Lambda E$ and there exists $e\in E$ with $t=r^\infty(e)$ and $\sigma(e)=e$.
\item Define $G=\Lambda E\,\sqcup\, \widehat G$.
\end{enumerate}
\fresult

Notice that the definition of $\sigma$ is in some sense `dual' to the definition of~$r$, because $r(t)$ is the join of all the irreducible elements which are strictly smaller than~$t$.
(However, the `true' dual of~$r$ is different~: it is the operation~$r$ in the opposite lattice~$T\op$, thus involving meet-irreducible elements.)
It is clear that the map $\sigma:T\to T$ is order-preserving.\par

In order to describe the effect of~$\sigma$, note first that $\sigma(t)=t$ if $t\in \Lambda E-E$ and $\sigma(t)>t$ if $t\notin\Lambda E$.
Now if $e\in E$, there are 3 cases~:
\begin{enumerate}
\item[(1)] If $]e,\hat1]\cap E$ has at least two minimal elements, then either $\sigma(e)=e$ or $\sigma(e)>e$ but $\sigma(e)$ is not irreducible.
\item[(2)]  If $]e,\hat1]\cap E$ has a unique minimal element~$e^+$, then $\sigma(e)=e^+$.
\item[(3)]  If $]e,\hat1]\cap E$ is empty (that is, $e$ is maximal in~$E$), then $\sigma(e)=\hat1$.
\end{enumerate}
Note that the equality $\sigma(e)=e$ also occurs in the third case for $e=\hat1$, provided $\hat1$ is irreducible.\par

Let $e\in E$ be as in the definition of~$t\in \widehat G$, that is, $t=r^\infty(e)$ and $\sigma(e)=e$.
Since the poset $]r^\infty(e),e]$ is totally ordered and consists of elements of~$E$,
it is clear that $r^\infty(f)=r^\infty(e)$ for every $f\in \; ]r^\infty(e),e]$ (see Lemma~\ref{interval}).
Since it may happen that $\sigma(f)=f$ for certain elements $f\in \; ]r^\infty(e),e]$,
the element $e$ in the definition above is not necessarily unique.
The following lemma shows that $e$ becomes unique if it is chosen minimal
among all elements $f\in E$ such that $r^\infty(f)=t$ and $\sigma(f)=f$.

\result{Lemma} \label{bulb} Let $t\in \widehat G$ and let $e\in E$ be minimal such that $t=r^\infty(e)$ and $\sigma(e)=e$.
Let $n\geq1$ be the smallest positive integer such that $r^n(e)=r^\infty(e)$,
so that $[t,e]= \{r^n(e),r^{n-1}(e), \ldots, r^1(e),e \}$.
\begin{enumerate}
\item $\sigma^i(t)=r^{n-i}(e)$ for $1\leq i \leq n$ and $\sigma^n(t)=\sigma^\infty(t)=e$.
\item $[t,\sigma^\infty(t)] =  \{r^n(e), \ldots, r(e), e \} =  \{t,\sigma(t), \ldots, \sigma^n(t) \}$.
\item $e=\sigma^\infty(t)$, in other words $e$ is unique. Moreover, $\sigma^\infty(t)\in E$.
\item $r^\infty(\sigma^i(t))=t$ for $1\leq i \leq n$.
\end{enumerate}
\fresult

\pf
By Lemma~\ref{interval}, $r^\infty(r^i(e)) = r^\infty(e)=t$ for all $0\leq i\leq n-1$.
If $1\leq i\leq n-1$, then $r^i(e)\in E$, $r^\infty(r^i(e)) =t$ and $r^i(e)<e$.
By minimality of~$e$, it follows that $\sigma(r^i(e))\neq r^i(e)$, hence $\sigma(r^i(e))> r^i(e)$.
Moreover, $\sigma(r^i(e))\leq r^{i-1}(e)$ by definition of~$\sigma(r^i(e))$ and the fact that $r^{i-1}(e)\in E$.
Since $[t,e]$ is totally ordered, this forces the equality $\sigma(r^i(e))= r^{i-1}(e)$.
This equality also holds if $i=n$ because $r^n(e)=t$ and $t\notin \Lambda E$, so $\sigma(t)>t$,
and again $\sigma(t)\leq r^{n-1}(e)$ so that $\sigma(t)= r^{n-1}(e)$.\par

Then one obtains $\sigma^i(t)=r^{n-i}(e)$ for $1\leq i \leq n$ and in particular $\sigma^n(t)=\sigma^\infty(t)=r^0(e)=e$.
The first three statements follow. The fourth is a consequence of Lemma~\ref{interval}.
\endpf

\result{Example} \label{Example-tot}
 {\rm If $T=\{0,1,\ldots,m\}$ is totally ordered, $E=\Irr(T)=\{1,\ldots,m\}$. Then $\widehat G=\{0\}$ and $G=T$.}
\fresult

\bigskip
We now show that the subset $G$ of~$T$ has another characterization.

\result{Lemma} \label{characterize-G} Let $G$ be as in Definition~\ref{notation-G}. Then
$$G=E\sqcup G^\sharp\mvirg$$
where $G^\sharp=\{a\in T\mid a=r^\infty \sigma^\infty(a)\}$.
\fresult

\pf First observe that $E\sqcup G^\sharp$ is a disjoint union because an element of the form $a=r^\infty \sigma^\infty(a)$ satisfies $r(a)=a$, so it cannot belong to~$E$.\par

In order to prove that $G\subseteq E\sqcup G^\sharp$, let $a\in G$. If $a\in E$, then obviously $a\in E\sqcup G^\sharp$. If $a\in\Lambda E-E$, then $a=\sigma(a)$, hence $a=\sigma^\infty(a)$. Moreover $a=r^\infty(a)$ since $a\notin E$. Hence $a=r^\infty \sigma^\infty(a)$, that is $a\in G^\sharp$. Finally if $a\in \widehat G$, then $a\notin E$ and $a=r^\infty \sigma^\infty(a)$, by Lemma~\ref{bulb}. This proves that $G\subseteq E\sqcup G^\sharp$.\par

For the reverse inclusion, first note that $E\subseteq G$ because $E\subseteq \Lambda E$. Now let $a\in G^\sharp$ and set $b=\sigma^\infty(a)$. If $b\notin E$, then $b=r^\infty(b)$, hence $b=a$ and $a=\sigma(a)$. It follows that $a\in \Lambda E$, hence $a\in G$. If now $b\in E$, there are two cases. Either $a\in\Lambda E$, hence $a\in G$ and we are done, or $a\notin \Lambda E$. But then we have $b=\sigma(b)$ and $a=r^\infty(b)$, so $a\in \widehat G$ by definition, hence $a\in G$. This proves the inclusion $E\sqcup G^\sharp\subseteq G$.
\endpf

The following two propositions will be crucial for our results on evaluations of fundamental functors in Sections~\ref{Section-generators} and~\ref{Section-independence}. We continue with the assumption that $T$ is a finite lattice and $(E,R)=\Irr(T)$.

\result{Proposition} \label{sigma-infty-in-E} Let $a\in T-\Lambda E$ and suppose that $\sigma^\infty(a)\in E$.
Let $m$ be the smallest positive integer such that $\sigma^\infty(a)=\sigma^m(a)$ and let $b:=r^\infty\sigma^\infty(a)$.
\begin{enumerate}
\item There exists $0\leq r\leq m-1$ such that $\sigma^r(a)<b<\sigma^{r+1}(a)$.
\item $b\in G$.
\item $]\sigma^r(a),\hat1]\cap E=[b,\hat1]\cap E=[\sigma^{r+1}(a),\hat1]\cap E$.
\end{enumerate}
\fresult

\pf
(a) Define $e_i=\sigma^i(a)$ for all $0\leq i\leq m$.
Note that $e_1,\ldots,e_{m-1}$ all belong to~$E$ because they belong to~$\Lambda E$
(since they are in the image of the operator~$\sigma$) and moreover $\sigma(e_i)> e_i$.
Also $e_m =\sigma^\infty(a)\in E$ by assumption.\par

We have $a=r^\infty(a) \leq r^\infty(\sigma^\infty(a)) < \sigma^\infty(a)$, because $\sigma^\infty(a) \in E$ by assumption.
Therefore, there is an integer $r\leq m-1$ such that $b\leq e_{r+1}$ but $b\not\leq e_r$.
The inequality $b< e_{r+1}$ is strict because $b\notin E$ while $e_{r+1}\in E$.
(The case $r=0$ occurs when $a\leq b< e_1$.)\par

In particular $b \leq r^\infty(e_{r+1}) \leq r^\infty(\sigma^\infty(a)) = b$, hence $b = r^\infty(e_{r+1})$.
Suppose that the element $e_r\vee b$ is irreducible.
Then either $e_r\vee b=b$ or $e_r\vee b=e_r$. The first case is impossible because
$b$ is not irreducible since $b= r^\infty(e_{r+1})\notin E$.
The second case is impossible because it would imply
$b \leq e_r$, contrary to the definition of~$r$.
Therefore $e_r\vee b\notin E$. Since $e_r\vee b \leq e_{r+1}$,
we obtain $e_r\vee b\leq r^\infty(e_{r+1})=b$ by definition of~$r^\infty(e_{r+1})$.
It follows that $e_r <b<e_{r+1}$, as required.\mpn

(b) To prove that $b\in G$, we first note that if $b\in \Lambda E$, then $b\in G$.
Otherwise, $b\notin \Lambda E$, $b=r^\infty(\sigma^\infty(a))$
and $\sigma(\sigma^\infty(a))=\sigma^\infty(a)$, proving that $b\in \widehat G$.\mpn

(c) By the definition of $\sigma(e_r)$, there is a unique minimal element in $]e_r,\hat1]\cap E$, namely $\sigma(e_r)=e_{r+1}$.
Therefore, $]e_r,\hat1]\cap E=[e_{r+1},\hat1]\cap E$. In particular $]e_r,\hat1]\cap E=[b,\hat1]\cap E$.
\endpf

\result{Notation} \label{def-zeta} 
Let $\zeta: G \longrightarrow \Iup(E,R)$ be the map defined by
$$\zeta(t)=\left\{ \begin{array}{ll} [t,\hat1] \,\cap\, E &\text{if }\; t\in E \,, \\
\rule{0ex}{2.5ex}]\sigma^\infty(t),\hat1]\cap E&\text{if }\; t\notin E \mpoint
\end{array}\right.$$
For any $B\in \Iup(E,R)$, define $\wedge B=\mathop{\wedge}_{e\in B}\limits e$.
By definition of~$\sigma^\infty(t)$, we obtain
$$\wedge\zeta(t)=\left\{ \begin{array}{ll} t &\text{if }\; t\in E \,, \\
\sigma^\infty(t) &\text{if }\; t\notin E\mpoint
\end{array}\right.$$
\fresult

\result{Proposition} \label{rho-leq} Let $t\in G$ and $t'\in T$ such that $t'\leq \wedge\zeta(t)$.
\begin{enumerate}
\item $\sigma^\infty(t')\leq \sigma^\infty(t)$ and $r^\infty(t')\leq r^\infty(t)$.
\item $t'\leq t$, except possibly if $t\in \widehat G$.
\item If $t'\not \leq t$, then $t\in \widehat G$ and $t=r^\infty(t')$.
\end{enumerate}
\fresult

\pf We have $G=(\Lambda E-E)\,\sqcup\, E \,\sqcup\, \widehat G$ and we consider the three cases for $t$ successively.\par

If $t\in \Lambda E-E$, then $\wedge\zeta(t)=\sigma^\infty(t)=t$, hence $t'\leq t$ and consequently $\sigma^\infty(t')\leq \sigma^\infty(t)$ and $r^\infty(t')\leq r^\infty(t)$.\par

If $t\in E$, then $\wedge\zeta(t)=t$, hence $t'\leq t$ and consequently $\sigma^\infty(t')\leq \sigma^\infty(t)$ and $r^\infty(t')\leq r^\infty(t)$.\par

Finally, if $t\in \widehat G$, then $\wedge\zeta(t)=\sigma^\infty(t)$, thus $\sigma^\infty(t')\leq \sigma^\infty(t)$. Moreover, using part~(d) of Lemma~\ref{rinfty} and part~(d) of Lemma~\ref{bulb}, we obtain
$$r^\infty(t')\leq r^\infty(\wedge\zeta(t))=r^\infty(\sigma^\infty(t))=t=r^\infty(t) \mpoint$$
This proves (a), and also (b) because the relation $t'\not \leq t$ can appear only if $t\in \widehat G$.
In that case case, by Lemma~\ref{bulb}, we can write $t=r^\infty(e)$ where $e=\sigma^\infty(t)$.
Since $t'\leq \sigma^\infty(t)$, we get
$t'\in \; [\hat0,t] \,\sqcup\, ]t,\sigma^\infty(t)]$ by Lemma~\ref{interval},
hence $t'\in \; ]t,\sigma^\infty(t)]$ because $t'\not \leq t$.
In other words, $t'\in\{\sigma(t), \ldots, \sigma^k(t) \}$ by Lemma~\ref{bulb}.
Therefore $t'= \sigma^i(t)$ for some $i\geq1$ and so $r^\infty(t')=t$ by Lemma~\ref{bulb}. This proves~(c) and completes the proof.
\endpf


\section{Correspondence functors} \label{Section-functors}

\medskip
\noindent
In this section, we recall the basic facts we need about correspondence functors. We refer to~\cite{BT2} for details.
We denote by $\CC$ the category of finite sets and correspondences.
Its objects are the finite sets and the set $\CC(Y,X)$ of morphisms from $X$ to~$Y$ is the set of all correspondences from $X$ to~$Y$,
namely all subsets of $Y\times X$ (using a reverse notation which is convenient for left actions).
If $S\subseteq Z\times Y$ and $R\subseteq Y\times X$, the composition of correspondences $SR$ is a correspondence from $X$ to~$Z$ defined by
$$SR=\{ (z,x)\in Z\times X \mid \exists \, y\in Y \,\text{ such that } \, (z,y)\in S \,\text{ and } \, (y,x)\in R \} \mpoint$$

When $X=Y$, a correspondence from $X$ to~$X$ is called a (binary) relation on~$X$, also called a Boolean matrix.
For any commutative ring~$k$, we let $k\CC$ be the $k$-linearization of~$\CC$.
The objects are the same, the set of morphisms $k\CC(Y,X)$ is the free $k$-module with basis $\CC(Y,X)$,
and composition is extended by $k$-bilinearity from composition in~$\CC$.
For any permutation $\sigma$ of~$X$, we write $\Delta_\sigma=\{ (\sigma(x),x)\mid x\in X \}$.
In particular, $\Delta_X:=\Delta_{\Id}$ is the identity morphism of the object~$X$.\par

A correspondence functor is a $k$-linear functor from $k\CC$ to the category $k\text{-\!}\Mod$ of left $k$-modules,
for some fixed commutative ring~$k$.
We let $\CF_k$ be the category of all correspondence functors.
If $F$ is a correspondence functor and $\psi\in k\CC(Y,X)$, we view the $k$-module homomorphism
$F(\psi): F(X)\to F(Y)$ as a left action of~$\psi$. More precisely, if $\alpha\in F(X)$, we define a left action $\psi\cdot\alpha:=F(\psi)(\alpha)\in F(Y)$.

If $E$ is a finite set, $\CC(E,E)$ is the monoid of all relations on~$E$ and we set
$$\CR_E:=k\CC(E,E) \mpoint$$
Among $\CR_E$-modules, there is the {\em fundamental module} $\CP_E f_R$,
associated to a poset $(E,R)$.
Here $\CP_E$ is a quotient algebra of the algebra $\CR_E$ and $f_R$ is a suitable idempotent in~$\CP_E$ depending on the order relation~$R$.
The algebra $\CP_E$ (called the algebra of permuted orders) and the module $\CP_E f_R$ were introduced in~\cite{BT1}.
All we need to know about $\CP_E f_R$ is its structure, described in the next result, which is Proposition~8.5 of~\cite{BT1} or Proposition~4.5 of~\cite{BT2}.

\result{Proposition} \label{fundamental-module}
Let $(E,R)$ be a finite poset.
\begin{enumerate}
\item The fundamental module $\CP_E f_R$ is a left module for the algebra $\CP_E$, hence also a left module for the algebra of relations~$\CR_E$.
\item $\CP_E f_R$ is a free $k$-module with a $k$-basis consisting of the elements $\Delta_\sigma f_R$,
where $\sigma$ runs through the group $\Sigma_E$ of all permutations of $E$.
\item $\CP_E f_R$ is a $(\CP_E,k\Aut(E,R))$-bimodule and the right action of $k\Aut(E,R)$ is free.
Explicitly, the right action of $\tau\in\Aut(E,R)$ maps the basis element $\Delta_\sigma f_R$ to the basis element $\Delta_{\sigma\tau} f_R$.
\item The action of the algebra of relations $\CR_E$ on the module $\CP_E f_R$
is given as follows. For any relation $Q\in \CC(E,E)$,
$$Q\cdot\Delta_\sigma f_R=\left\{\begin{array}{ll}
\Delta_{\tau\sigma}f_R&\hbox{if}\;\;\exists\tau\in\Sigma_E\;\hbox{such that}\;
\Delta_E\subseteq \Delta_{\tau^{-1}}Q\subseteq {\ls\sigma R},\\
0&\hbox{otherwise}\mvirg\end{array}\right.$$
where $\ls\sigma R=\big\{\big(\sigma(e),\sigma(f)\big)\mid (e,f)\in R\big\}$, or equivalently $\ls\sigma R=\Delta_\sigma R\Delta_{\sigma^{-1}}$.
(Note that $\tau$ is unique in the first case.)
\end{enumerate}
\fresult

The evaluation at~$E$ of a correspondence functor is a left $\CR_E$-module and our strategy will be to work with correspondence functors rather than $\CR_E$-modules.
In particular, given a finite poset $(E,R)$, we defined in~\cite{BT2} a {\em fundamental functor} $\S_{E,R}$ which is associated with the fundamental module $\CP_E f_R$ and has the following properties (see Proposition~2.6 in~\cite{BT3}).

\result{Proposition} \label{SER} Let $(E,R)$ be a finite poset and $X$ a finite set.
\begin{enumerate}
\item $\S_{E,R}(X)=\{0\}$ if $|X|<|E|$.
\item $\S_{E,R}(E)\cong\CP_E f_R$.
\end{enumerate}
\fresult

For the largest part of the present paper, we do not need to go back to the definition of the fundamental functor $\S_{E,R}$, because it is fully described by Theorem~\ref{kertheta} below.
However, the precise construction of $\S_{E,R}$ will be recalled in Section~\ref{Section-simple-tensor}, where we will need to analyze a right $k\Aut(E,R)$-module structure on~$\S_{E,R}(X)$ (induced by the right $k\Aut(E,R)$-module structure on~$\CP_E f_R$).

Now we explain the connection between correspondence functors and lattices.
Let $T$ be a finite lattice. We defined in~\cite{BT3} a correspondence functor~$F_T$ as follows.
If $X$ is a finite set, then $F_T(X)=kT^X$, the free $k$-module with basis the set $T^X$ of all functions from $X$ to~$T$.
If $R\subseteq Y\times X$ is a correspondence and if $\varphi\in T^X$, then 
we associate the function $R\cdot\varphi=F_T(R)(\varphi) \in T^Y$, also simply written $R\varphi$, defined by
$$(R\varphi)(y):=\mbigvee{(y,x)\in R}\varphi(x)\mvirg$$
with the usual rule that a join over the empty set is equal to~$\hat{0}$. The map
$$F_T(R):F_T(X)\longrightarrow F_T(Y)$$
is the unique $k$-linear extension of this construction.\par

If $E$ denotes the set of irreducible elements in~$T$, the functions $f\in F_T(X)$ such that $E\not\subseteq f(X)$ generate a subfunctor, written~$H_T$ in~\cite{BT3}, which we often use.
The functors $F_T$ play an important role because they are connected to fundamental functors by a morphism described in the following result (see Theorem~6.5 in~\cite{BT3}).

\result{Theorem} \label{surjection} Let $T$ be a finite lattice, let $(E,R)=\Irr(T)$, and let $\iota:E\to T$ denote the inclusion map.
\begin{enumerate}
\item There exists a unique surjective morphism of correspondence functors
$$\Theta_T:F_T\longrightarrow \S_{E,R\op}$$
such that $\Theta_{T,E}(\iota) = f_{R\op}$ (an element in $\S_{E,R\op}(E)=\CP_E f_{R\op}$).
\item The subfunctor $H_T$ is contained in the subfunctor $\Ker(\Theta_T)$.
Explicitly, if $X$ is a finite set and if $f\in F_T(X)$ satisfies the condition $E\not\subseteq f(X)$, then $\Theta_{T,X}(f)=0$.
\item The functor $F_T$ is generated by $\iota\in F_T(E)$, while the functor $\S_{E,R\op}$ is generated by $f_{R\op}\in\S_{E,R\op}(E)$.
\end{enumerate}
\fresult

In order to have control of the fundamental functor $\S_{E,R\op}$, we need to understand the kernel of~$\Theta_T$.
To this end, we need to consider some correspondences which were introduced in~\cite{BT3} and which play again an important role in the present paper.

\result{Notation} \label{Gamma}
Let $T$ be a finite lattice and $(E,R)=\Irr(T)$.
For any finite set $X$ and any map $\varphi: X\to T$, we associate the correspondence
$$\Gamma_\varphi:= \{ (x,e)\in X\times E \mid e\leq_T \varphi(x) \} \subseteq X\times E \mpoint$$
In the special case where $T=\Idown(E,R)$, we obtain
$$\Gamma_\varphi= \{ (x,e)\in X\times E \mid e\in \varphi(x) \} \mpoint$$
\fresult

For the description of the kernel of~$\Theta_T:F_T\to\S_{E,R\op}$,
the following result was obtained as Theorem~7.1 in~\cite{BT3}.
The result actually gives a first explicit description of every fundamental functor and it is one of our main tools in this paper.

\result{Theorem} \label{kertheta} Let $T$ be a finite lattice, let $(E,R)=\Irr(T)$, and let $X$ be a finite set.
The kernel of the map
$$\Theta_{T,X}:F_T(X)\longrightarrow \S_{E,R\op}(X)$$
is equal to the set of linear combinations $\sum_{\varphi:X\to T}\limits\lambda_\varphi\varphi$, where $\lambda_\varphi\in k$, such that for any map $\psi:X\to \Iup(E,R)$
$$\sumb{\varphi}{\Gamma_\psi\op\Gamma_\varphi = R\op}\lambda_\varphi=0\mpoint$$
Here $\Gamma_{\varphi}=\{(x,e)\in X\times E\mid e \leq_T \varphi(x)\}$
and $\Gamma_\psi\op=\big\{(e,x)\in E\times X\mid e\in\psi(x)\big\}$, as in Notation~\ref{Gamma}.
\fresult

In order to use the condition $\Gamma_\psi\op\Gamma_\varphi = R\op$ appearing in Theorem~\ref{kertheta}, we shall also need equivalent formulations.
We first fix notation. If $\psi:X\to \Iup(E,R)$ is a map, define the function $\wedge\psi:X\to T$ by
$$\forall x\in X,\;\;\wedge\psi(x)=\mbigwedge{e\in\psi(x)}e \mvirg$$
where $\bigwedge$ is the meet in the lattice~$T$. If $\varphi$ and $\varphi'$ are two functions $X\to T$, we write $\varphi\leq \varphi'$ if $\varphi(x)\leq \varphi'(x)$ for all $x\in X$.
The following result is Theorem~7.3 in~\cite{BT3}. 

\result{Theorem} \label{Gammapsigamma}
Let $T$ be a finite lattice, let $(E,R)=\Irr(T)$, let $\iota:E\to T$ denote the inclusion map, and let $X$ be a finite set.
Let $\varphi:X \to T$ and $\psi:X\to \Iup(E,R)$ be maps with associated correspondence
$\Gamma_\varphi$ and $\Gamma_\psi\op$, as in Theorem~\ref{kertheta} above.
The following conditions are equivalent.
\begin{enumerate}
\item $\Gamma_\psi\op\varphi =\iota$.
\item $\Gamma_\psi\op\Gamma_\varphi \iota=\iota$.
\item $\Delta_E\subseteq\Gamma_\psi\op\Gamma_\varphi\subseteq R\op$.
\item $\Gamma_\psi\op\Gamma_\varphi= R\op$.
\item $\varphi\leq\wedge\psi$ and
$\forall e\in E,\;\exists x\in X$ such that $\varphi(x)=e$ and $\psi(x)=[e,\cdot[_E$.
\item $\forall t\in T,\;\psi\big(\varphi^{-1}(t)\big)\subseteq[t,\cdot[_T\cap E$ and
$\;\forall e\in E,\;\psi\big(\varphi^{-1}(e)\big)=[e,\cdot[_E$.
\end{enumerate}
\fresult

Conditions (e) and (f) will play a crucial role in our results on fundamental functors and simple functors (Sections~\ref{Section-generators} and~\ref{Section-independence}).

The next step in our strategy is to realize every simple correspondence functor as a quotient of a fundamental functor.
More precisely, given a finite poset $(E,R)$ and a left $k\Aut(E,R)$-module~$V$,
we defined in~\cite{BT2} a correspondence functor $S_{E,R,V}$ with the following properties (see Proposition~2.6 and Lemma~2.7 in~\cite{BT3}).

\result{Proposition} \label{SERV} Let $(E,R)$ be a finite poset, let $V$ be a left $k\Aut(E,R)$-module generated by a single element~$v$, and let $X$ be a finite set.
\begin{enumerate}
\item $S_{E,R,V}(X)=\{0\}$ if $|X|<|E|$.
\item $S_{E,R,V}(E)=\CP_E f_R\otimes_{k\Aut(E,R)}V$.
\item There is a surjective morphism of correspondence functors
$$\Phi: \S_{E,R} \longrightarrow S_{E,R,V}$$
such that, on evaluation at the finite set~$E$,
we obtain the surjective homomorphism of $\CR_E$-modules
$$\Phi_E: \CP_E f_R \longrightarrow \CP_E f_R\otimes_{k\Aut(E,R)}V \,, \qquad a\mapsto a\otimes v \mpoint$$
\end{enumerate}
\fresult

The precise construction of $S_{E,R,V}$ will be recalled in Section~\ref{Section-simple-tensor} and the definition of~$\Phi$ will be given in Notation~\ref{Notation-Psi}.
In the special case when $k$ is a field and the $k\Aut(E,R)$-module $V$ is simple, we obtain more (see Theorem~4.7 in~\cite{BT2}).

\result{Theorem} \label{parametrization} 
Assume that $k$ is a field.
\begin{enumerate}
\item If the $k\Aut(E,R)$-module $V$ is simple, then the functor $S_{E,R,V}$ is simple and $S_{E,R,V}(E)=\CP_E f_R\otimes_{k\Aut(E,R)}V$ is a simple $\CR_E$-module.
\item The map $(E,R,V)\mapsto S_{E,R,V}$ provides a parametrization of all simple functors by isomorphism classes of triples $(E,R,V)$, where $(E,R)$ is a finite poset and $V$ is a simple $k\Aut(E,R)$-module.
\end{enumerate}
\fresult

In order to obtain information about simple functors $S_{E,R,V}$, we shall always work first with the fundamental functor $\S_{E,R}$, which is a precursor of~$S_{E,R,V}$ since we recover $S_{E,R,V}$ by means of the surjective morphism $\Phi: \S_{E,R} \to S_{E,R,V}$.
This explains why the fundamental functors play a crucial role throughout our work.
We shall see in Section~\ref{Section-simple-tensor} that there is an explicit way to recover $S_{E,R,V}$ from~$\S_{E,R}$.
It is also worth mentioning that both $\S_{E,R}$ and $S_{E,R,V}$ are defined over an arbitrary commutative ring~$k$.


\section{Generators for the evaluations of fundamental functors} \label{Section-generators}

\medskip
\noindent
As usual, $E$ denotes a fixed finite set and $R$ an order relation on~$E$.
Our purpose is to prove that, for an arbitrary commutative ring~$k$ and for any finite set $X$, the evaluation $\S_{E,R}(X)$ of the fundamental correspondence functor $\S_{E,R}$ is a free $k$-module, by finding an explicit $k$-basis. In this section, we first deal with $k$-linear generators.\par

Let $T$ be any lattice such that $(E,R)=\Irr(T)$.
Note that $\Idown(E,R)$ is the largest such lattice and that any other is a quotient of~$\Idown(E,R)$ (Lemma~\ref{principal}).
By Theorem~\ref{surjection}, the fundamental functor $\S_{E,R\op}$ is isomorphic to a quotient of~$F_T$ via a morphism
$$\Theta_T:F_T \longrightarrow \S_{E,R\op} \mpoint$$
For this reason, we work with $\S_{E,R\op}$ rather than~$\S_{E,R}$.

\result{Notation} \label{notation-BX}
Let $G=G(T)$ be the subset defined in Notation~\ref{notation-G} and let $X$ be a finite set. We define
$\CB_X$ to be the set of all maps $\varphi:X \to T$ such that $E\subseteq \varphi(X) \subseteq G$.
\fresult

Our main purpose is to prove that the set $\Theta_{T,X}(\CB_X)$ is a $k$-basis of~$\S_{E,R\op}(X)$.
We first prove in this section that $\Theta_{T,X}(\CB_X)$ generates $k$-linearly $\S_{E,R\op}(X)$ and
then we shall show in Section~\ref{Section-independence} that $\Theta_{T,X}(\CB_X)$ is $k$-linearly independent.\par

By Lemma~\ref{characterize-G}, we have $G=E\sqcup G^\sharp$ where $G^\sharp=\{a\in T \mid a=r^\infty\sigma^\infty(a)\}$,
with $r^\infty$ defined in Notation~\ref{r(t)} and $\sigma^\infty$ in Definition~\ref{notation-G}.
We denote by $G^c$ the complement of $G$ in~$T$, namely
$$G^c=\{a\in T\mid a\notin E,\;a<r^\infty \sigma^\infty(a)\}\mpoint$$

\result{Lemma} \label{reduction sequence}
Let $a\in G^c$, and let $b=r^\infty \sigma^\infty(a)$. There exists an integer $r\geq 0$ such that
$$a<\sigma(a)<\ldots<\sigma^r(a)<b \leq \sigma^{r+1}(a)\mvirg$$
and $\sigma^j(a)\in E$ for $j\in\{1,\ldots,r\}$.
Moreover $b=r^\infty \sigma^\infty\big(\sigma^j(a)\big)$ for $j\in\{1,\ldots,r\}$, and $b\in G^\sharp$, where $G=E\sqcup G^\sharp$ as in Lemma~\ref{characterize-G}.
\fresult

\pf
We know that $a=\sigma^0(a)\notin E$ because $a\in G^c$.
Also $a<\sigma(a)$ because $a\notin \Lambda E \subseteq G$.
Suppose first that there exists an integer $r\geq 0$ such that $c=\sigma^{r+1}(a)\notin E$.
In this case, we choose $r$ minimal with this property, so that $\sigma^j(a)\in E$ for $1\leq j\leq r$.
We have $c=\sigma\big(\sigma^{r}(a)\big)\in\Lambda E-E$, hence $c=\sigma(c)=\sigma^\infty(c)=\sigma^\infty(a)$.
Moreover $b=r^\infty(c)=c$, because $c\notin E$.
Since $\sigma(c)=c$, we obtain $\sigma(b)=b$, hence $b=r^\infty \sigma^\infty(b)$.
Therefore $b\in G^\sharp$.\par

Suppose now that $\sigma^r(a)\in E$ for all $r\in\Z_{>0}$. Then Proposition~\ref{sigma-infty-in-E} applies and
there exists an integer $r\geq 0$ such that $\sigma^r(a)<b<\sigma^{r+1}(a)$ and $b\in G$.
Moreover, $b\notin E$, because $b=r^\infty(b)$, so $b\in G^\sharp$.
\endpf

\result{Definition} \label{def reduction sequence}
For $a\in G^c$, the sequence $a<\sigma(a)<\ldots<\sigma^r(a)<b$ defined in Lemma~\ref{reduction sequence} will be called the {\em reduction sequence} associated to~$a$. 
\fresult

\result{Notation} \label{sequences and maps} Let $n\geq 1$ and let $(a_0,a_1,\ldots,a_n)$ be a sequence of distinct elements of~$T$.
We denote by $[a_0,\ldots,a_n]:T\to T$ the map defined by
$$\forall t\in T,\;\;[a_0,\ldots,a_n](t)=\left\{\begin{array}{ll}a_{j+1}&\hbox{if}\;t=a_j,\;j\in\{0,\ldots,n-1\}\\t&\hbox{otherwise.}\end{array}\right.$$
If $a\in G^c$, let $(a_0,a_1,\ldots,a_r,a_{r+1})$ be the reduction sequence associated to~$a$, with $a_0=a$ and $a_{r+1}=b=r^\infty \sigma^\infty(a)$.
We then denote by $u_a$ the element of $k(T^T)=F_T(T)$ defined by
$$u_a=[a_0,a_1]-[a_0,a_1,a_2]+\ldots+(-1)^r[a_0,a_1,\ldots,a_{r+1}]\mpoint$$

\fresult

We can now describe a family of useful elements in $\Ker\Theta_T$.

\result{Theorem} \label{phi minus uaphi}
Let $T$ be a finite lattice, let $(E,R)=\Irr(T)$, let $G^c=\{a\in T\mid a\notin E, \; a<r^\infty \sigma^\infty(a)\}$, and let $X$ be a finite set.
Then, for any $a\in G^c$ and for any function $\varphi:X\to T$,
$$\varphi - u_a\circ\varphi\in\Ker\Theta_{T,X}\mvirg$$
where $u_a\circ\varphi$ is defined by bilinearity from the composition of maps $T^T\times T^X\to T^X$.
\fresult

\pf The kernel of the map $\Theta_{T,X}:F_T(X)\to \S_{E,R\op}(X)$ was described in Theorem~\ref{kertheta}. Let $\sum_{\varphi:T\to X}\limits\lambda_\varphi \varphi\in F_T(X)$, where $\lambda_\varphi\in k$. Then $\sum_{\varphi:T\to X}\limits\lambda_\varphi \varphi\in\Ker\Theta_{T,X}$ if and only if the coefficients $\lambda_\varphi$ satisfy a system of linear equations indexed by maps $\psi:X\to \Iup(E,R)$.
The equation $(E_\psi)$ indexed by such a map $\psi$ is the following~:
$$(E_\psi): \qquad \sum_{\varphi\vdashER{E,R}\psi}\lambda_\varphi=0\mvirg$$
where $\varphi\vdashER{E,R}\psi$ means that $\varphi:X\to T$ and $\psi:X\to \Iup(E,R)$ satisfy the equivalent conditions of Theorem~\ref{Gammapsigamma}. We shall use condition~(f) of Theorem~\ref{Gammapsigamma}, namely
$$
\varphi\vdashER{E,R}\psi\iff 
\left\{\begin{array}{l}\forall t\in T,\;\;\psi\big(\varphi^{-1}(t)\big)\subseteq[t,\cdot[_T\cap E \mvirg\\
\rule{0ex}{3ex}\forall e\in E,\;\;\psi\big(\varphi^{-1}(e)\big)=[e,\cdot[_E\mpoint\end{array}\right.
$$
Let $a\in  G^c$, and let $(a,e_1,e_2,\ldots,e_r,b)$ be the associated reduction sequence.
Recall that $e_1,\ldots,e_r\in E$ but $b\notin E$.
If $r\geq 1$, note that $[a,\cdot[_T\cap E=[e_1,\cdot[_E$ because $a<\sigma(a)=e_1\in E$.
Define, for each $i\in\{1,\ldots ,r\}$,
$$\varphi_i=[a,e_1,\ldots,e_i]\circ\varphi \,, \qquad\text{and also } \;\varphi_{r+1}=[a,e_1,\ldots,e_r,b]\circ\varphi \mpoint$$
In particular, for any $i\in\{1,\ldots,r+1\}$,
$$\text{if } \varphi(x) \in T-\{a,e_1,\ldots,e_r\} \mvirg \;\text{ then } \varphi_i(x)=\varphi(x) \mpoint$$
The other values of the maps $\varphi_i$ are given in the following table~:
\begin{equation}\label{values of phi}
\begin{array}{c||c|c|c|c|c|c|}
x\in &\varphi^{-1}(a)&\varphi^{-1}(e_1)&\varphi^{-1}(e_2)&\ldots&\varphi^{-1}(e_{r-1})&\varphi^{-1}(e_r)\\
\hline
\hline
\varphi(x)&a&e_1&e_2&\ldots&e_{r-1}&e_r\\
\hline
\varphi_1(x)&e_1&e_1&e_2&\ldots&e_{r-1}&e_r\\
\hline
\varphi_2(x)&e_1&e_2&e_2&\ldots&e_{r-1}&e_r\\
\hline
\varphi_3(x)&e_1&e_2&e_3&\ldots&e_{r-1}&e_r\\
\hline
\ldots&\ldots&\ldots&\ldots&\ldots&\ldots&\ldots\\
\hline
\varphi_r(x)&e_1&e_2&e_3&\ldots&e_{r}&e_r\\
\hline
\varphi_{r+1}(x)&e_1&e_2&e_3&\ldots&e_{r}&b\\
\hline
\end{array}
\end{equation}

\bigskip
\noindent
We want to prove that the element
$$\varphi-u_a\circ\varphi=\varphi-\varphi_1+\varphi_2-\ldots +(-1)^{r-1}\varphi_{r+1}$$
belongs to $\Ker\Theta_{T,X}$. We must prove that it satisfies the equation $(E_\psi)$ for every~$\psi$,
so we must find which of the functions $\varphi,\varphi_1,\varphi_2, \ldots, \varphi_{r+1}$
are linked with~$\psi$ under the relation~$\vdashER{E,R}$.
We are going to prove that only two consecutive functions can be linked with a given~$\psi$, from which it follows that the corresponding equation $(E_\psi)$ is satisfied because it reduces to either $1-1=0$, or $-1+1=0$.
Of course, if none of $\varphi,\varphi_1,\varphi_2, \ldots, \varphi_{r+1}$ is linked with~$\psi$, then the corresponding equation $(E_\psi)$ is just $0=0$.
It follows from this that $\varphi-u_a\circ\varphi$ satisfies all equations $(E_\psi)$, hence belongs to $\Ker\Theta_{T,X}$, as required.
We note for completeness that it may happen that some of the functions $\varphi,\varphi_1,\varphi_2, \ldots, \varphi_{r+1}$ are equal (this occurs if an inverse image is empty in some column of the table), but this does not play any role in the argument.\par

Assume first that $r\geq1$.
Write $U:=T-\{a,e_1,\ldots,e_r\}$ and $V:=E-\{e_1,\ldots,e_r\}$.
The linking with a fixed~$\psi$ is controlled by the following conditions~:
$$\varphi\vdashER{E,R} \psi \iff \left\{\begin{array}{l}\forall t\in U,\;\;\psi\big(\varphi^{-1}(t)\big)\subseteq [t,\cdot[_T\cap E\\\rule{0ex}{3ex}\forall e\in V,\;\;\psi\big(\varphi^{-1}(e)\big)=[e,\cdot[_E\\\rule{0ex}{3ex}\psi\big(\varphi^{-1}(a)\big)\subseteq [a,\cdot[_T\cap E=[e_1,\cdot[_E\\
\rule{0ex}{3ex}\psi\big(\varphi^{-1}(e_i)\big)=[e_i,\cdot[_E\;\;\forall i\in\{1,\ldots,r\}\mpoint\end{array}\right.$$
The subsets $\varphi_j^{-1}(e_i)$ are determined by Table~\ref{values of phi} and can be written in terms of~$\varphi$.
In particular, $\varphi_1^{-1}(e_1)=\varphi^{-1}(a)\sqcup\varphi^{-1}(e_1)$, so we obtain
$$\varphi_1\vdashER{E,R} \psi \iff \left\{\begin{array}{l}\forall t\in U,\;\;\psi\big(\varphi^{-1}(t)\big)\subseteq [t,\cdot[_T\cap E\\\rule{0ex}{3ex}\forall e\in V,\;\;\psi\big(\varphi^{-1}(e)\big)=[e,\cdot[_E\\\rule{0ex}{3ex}\psi\big(\varphi^{-1}(a)\sqcup\varphi^{-1}(e_1)\big)= [e_1,\cdot[_E\\\rule{0ex}{3ex}\psi\big(\varphi^{-1}(e_i)\big)=[e_i,\cdot[_E\;\;\forall i\in\{2,\ldots,r\}\mpoint\end{array}\right.$$
Similarly, for $2\leq j \leq r$, we have $\varphi_j^{-1}(e_1)=\varphi^{-1}(a)$ and $\varphi_j^{-1}(e_{i+1})=\varphi^{-1}(e_i)$ if $1\leq i\leq j-2$, and then
$\varphi_j^{-1}(e_j)=\varphi^{-1}(e_{j-1})\sqcup\varphi^{-1}(e_j)$.
Therefore we get successively
$$\varphi_2\vdashER{E,R} \psi \iff \left\{\begin{array}{l}\forall t\in U,\;\;\psi\big(\varphi^{-1}(t)\big)\subseteq [t,\cdot[_T\cap E\\
\rule{0ex}{3ex}\forall e\in V,\;\;\psi\big(\varphi^{-1}(e)\big)=[e,\cdot[_E\\
\rule{0ex}{3ex}\psi\big(\varphi^{-1}(a)\big)=[e_1,\cdot[_E\\
\rule{0ex}{3ex}\psi\big(\varphi^{-1}(e_1)\sqcup\varphi^{-1}(e_2)\big)= [e_2,\cdot[_E\\
\rule{0ex}{3ex}\psi\big(\varphi^{-1}(e_i)\big)=[e_i,\cdot[_E\;\;\forall i\in\{3,\ldots,r\}\mpoint
\end{array}\right.$$
\vspace{1ex}
$$\varphi_3\vdashER{E,R} \psi \iff \left\{\begin{array}{l}\forall t\in U,\;\;\psi\big(\varphi^{-1}(t)\big)\subseteq [t,\cdot[_T\cap E\\
\rule{0ex}{3ex}\forall e\in V,\;\;\psi\big(\varphi^{-1}(e)\big)=[e,\cdot[_E\\
\rule{0ex}{3ex}\psi\big(\varphi^{-1}(a)\big)=[e_1,\cdot[_E\\
\rule{0ex}{3ex}\psi\big(\varphi^{-1}(e_1)\big)=[e_2,\cdot[_E\\
\rule{0ex}{3ex}\psi\big(\varphi^{-1}(e_2)\sqcup\varphi^{-1}(e_3)\big)= [e_3,\cdot[_E\\
\rule{0ex}{3ex}\psi\big(\varphi^{-1}(e_i)\big)=[e_i,\cdot[_E\;\;\forall i\in\{4,\ldots,r\}\mpoint
\end{array}\right.$$
\vspace{1ex}
$$\ldots$$
\vspace{1ex}
$$\varphi_{r-1}\vdashER{E,R} \psi \iff \left\{\begin{array}{l}\forall t\in U,\;\;\psi\big(\varphi^{-1}(t)\big)\subseteq [t,\cdot[_T\cap E\\
\rule{0ex}{3ex}\forall e\in V,\;\;\psi\big(\varphi^{-1}(e)\big)=[e,\cdot[_E\\
\rule{0ex}{3ex}\psi\big(\varphi^{-1}(a)\big)=[e_1,\cdot[_E\\
\rule{0ex}{3ex}\psi\big(\varphi^{-1}(e_i)\big)=[e_{i+1},\cdot[_E\;\forall i\in\{1,\ldots, r-3\}\\
\rule{0ex}{3ex}\psi\big(\varphi^{-1}(e_{r-2})\sqcup\varphi^{-1}(e_{r-1})\big)= [e_{r-1},\cdot[_E\\
\rule{0ex}{3ex}\psi\big(\varphi^{-1}(e_r)\big)=[e_r,\cdot[_E\mpoint
\end{array}\right.$$
\vspace{1ex}
$$\varphi_r\vdashER{E,R} \psi \iff \left\{\begin{array}{l}\forall t\in U,\;\;\psi\big(\varphi^{-1}(t)\big)\subseteq [t,\cdot[_T\cap E\\
\rule{0ex}{3ex}\forall e\in V,\;\;\psi\big(\varphi^{-1}(e)\big)=[e,\cdot[_E\\
\rule{0ex}{3ex}\psi\big(\varphi^{-1}(a)\big)=[e_1,\cdot[_E\\
\rule{0ex}{3ex}\psi\big(\varphi^{-1}(e_i)\big)=[e_{i+1},\cdot[_E\;\;\forall i\in\{1,\ldots,r-2\}\\
\rule{0ex}{3ex}\psi\big(\varphi^{-1}(e_{r-1})\sqcup\varphi^{-1}(e_r)\big)= [e_r,\cdot[_E\mpoint\\
\end{array}\right.$$
\vspace{1ex}
$$\varphi_{r+1}\vdashER{E,R} \psi \iff \left\{\begin{array}{l}\forall t\in U,\;\;\psi\big(\varphi^{-1}(t)\big)\subseteq [t,\cdot[_T\cap E\\
\rule{0ex}{3ex}\forall e\in V,\;\;\psi\big(\varphi^{-1}(e)\big)=[e,\cdot[_E\\
\rule{0ex}{3ex}\psi\big(\varphi^{-1}(a)\big)=[e_1,\cdot[_E\\
\rule{0ex}{3ex}\psi\big(\varphi^{-1}(e_i)\big)=[e_{i+1},\cdot[_E\;\;\forall i\in\{1,\ldots,r-1\}\\
\rule{0ex}{3ex}\psi\big(\varphi^{-1}(e_r)\big)\subseteq[b,\cdot[_T\cap E\mpoint
\end{array}\right.$$\vspace{1ex}
\par

Suppose that $\varphi\vdashER{E,R} \psi$. This clearly implies $\varphi_1\vdashER{E,R} \psi$. Also $\varphi_i\nvdashER{E,R} \psi$ for $i\geq 2$, because $\varphi\vdashER{E,R} \psi$ implies $\psi\big(\varphi^{-1}(e_1)\big)=[e_1,\cdot[_E$, but $\varphi_i\vdashER{E,R} \psi$ implies $\psi\big(\varphi^{-1}(e_1)\big)\subseteq [e_2,\cdot[_E$ when $i\geq 2$. Therefore only $\varphi$ and $\varphi_1$ are involved in this case.\par

Suppose now that $\varphi_1\vdashER{E,R} \psi$ but $\varphi\nvdashER{E,R} \psi$. Then $\psi\big(\varphi^{-1}(e_1)\big)\subseteq\, ]e_1,\cdot[_T\cap E=[e_2,\cdot[_E$ (because {$e_2=\sigma(e_1)$}) and $\psi\big(\varphi^{-1}(a)\big)=[e_1,\cdot[_E$, hence in particular $\varphi_2\vdashER{E,R} \psi$, since $\varphi_1\vdashER{E,R} \psi$ also implies {$\psi\big(\varphi^{-1}(e_i)\big)=[e_i,\cdot[_E$} for $i\in\{2,\ldots,n\}$.
On the other hand, since $\varphi_1\vdashER{E,R} \psi$ implies $\psi\big(\varphi^{-1}(e_2)\big)=[e_2,\cdot[_E$,
we cannot have $\psi\big(\varphi^{-1}(e_2)\big)\subseteq [e_3,\cdot[_E$ and so
$\varphi_i\nvdashER{E,R} \psi$, for $i\geq 3$. Therefore only $\varphi_1$ and $\varphi_2$ are involved in this case.\par

Suppose by induction that $\varphi_i\vdashER{E,R} \psi$ but $\varphi_{i-1}\nvdashER{E,R} \psi$, for some $i\in \{1,\ldots,r-1\}$. Then the same argument shows that $\varphi_{i+1}\vdashER{E,R} \psi$ and that only $\varphi_i$ and $\varphi_{i+1}$ are involved in this case.\par

Suppose now that $\varphi_r\vdashER{E,R} \psi$ but $\varphi_{r-1}\nvdashER{E,R}\psi$. Then $\psi\big(\varphi^{-1}(e_{r-1})\sqcup\varphi^{-1}(e_r)\big)= [e_r,\cdot[_E$ but $\psi\big(\varphi^{-1}(e_r)\big)\neq [e_r,\cdot[_E$. Hence $\psi\big(\varphi^{-1}(e_r)\big)\subseteq\, ]e_r,\cdot[_E\,\subseteq [b,\cdot[_T\cap E$, since $\sigma(e_r)\geq b$. Moreover $e_r\in\psi\big(\varphi^{-1}(e_{r-1})\sqcup\varphi^{-1}(e_r)\big)$ and $e_r\notin\psi\big(\varphi^{-1}(e_{r})\big)$. It follows that $e_r\in\psi\big(\varphi^{-1}(e_{r-1})\big)$, hence $\psi\big(\varphi^{-1}(e_{r-1})\big)=[e_r,\cdot[_E$.
Therefore $\varphi_{r+1}\vdashER{E,R} \psi$ and only $\varphi_r$ and $\varphi_{r+1}$ are involved in this case.\par

Finally, if $\varphi_{r+1}\vdashER{E,R} \psi$, then $\psi\big(\varphi^{-1}(e_{r-1})\big)=[e_r,\cdot[_E$ and $\psi\big(\varphi^{-1}(e_{r})\big)\subseteq [b,\cdot[_T\cap E\subseteq [e_r,\cdot[_E$. Thus $\psi\big(\varphi^{-1}(e_{r-1})\sqcup \varphi^{-1}(e_{r})\big)=[e_r,\cdot[_E$.
Therefore $\varphi_r\vdashER{E,R} \psi$ and we are again in the case when only $\varphi_r$ and $\varphi_{r+1}$ are involved.\par

The special case $r=0$ has to be treated separately. There are only 2 terms $\varphi$ and $\varphi_1$ in the alternating sum.
If $\varphi\vdashER{E,R} \psi$, then $\psi\big(\varphi^{-1}(a)\big)\subseteq [a,\cdot[_T\cap E$,
hence $\psi\big(\varphi^{-1}(a)\big)\subseteq [b,\cdot[_T\cap E$ because $b\leq\sigma(a)$ when $r=0$. Therefore
$$\psi\big(\varphi_1^{-1}(b)\big)=\psi\big(\varphi^{-1}(a)\sqcup \varphi^{-1}(b)\big) \subseteq [b,\cdot[_T\cap E$$
and so $\varphi_1\vdashER{E,R} \psi$.
Conversely, it is straightforward to see that $\varphi_1\vdashER{E,R} \psi$ implies $\varphi\vdashER{E,R} \psi$.\par

We have proved that only two consecutive functions can be linked with a given~$\psi$, as was to be shown.
\endpf

We have now paved the way for finding generators of $\S_{E,R\op}(X)$.

\result{Theorem} \label{generators-SER} Let $T$ be a finite lattice, let $(E,R)=\Irr(T)$, and let
$$G^c=\{a\in T\mid  a\notin E, \; a<r^\infty \sigma^\infty(a)\} \mpoint$$
For $a\in G^c$, let $u_a$ be the element of~$k(T^T)$ introduced in Notation~\ref{sequences and maps}, and let $u_T$ denote the composition of all the elements $u_a$, for $a\in G^c$, in some order (they actually commute, see Theorem~\ref{idempotents-ua} below).
\begin{enumerate}
\item Let $X$ be a finite set. Then for any $\varphi:X\to T$, the element $u_T\circ\varphi$ is a $k$-linear combination of functions $f:X\to T$ such that $f(X)\subseteq G$.
\item Let $\CB_X$ be the set of all maps $\varphi:X \to T$ such that $E\subseteq \varphi(X) \subseteq G$.
Then the set $\Theta_{T,X}(\CB_X)$ generates $\S_{E,R\op}(X)$ as a $k$-module.
\end{enumerate}
\fresult

\pf (a) We see in Table~\ref{values of phi} that the functions $\varphi_1,\varphi_2,\ldots,\varphi_{r+1}$ do not take the value~$a$.
It follows that for any $a\in G^c$ and any $\varphi:X\to T$, the element $u_a\circ \varphi=\varphi_1-\varphi_2+\ldots+(-1)^{r}\varphi_{r+1}$ is a $k$-linear combination of functions $\varphi_i$ such that $\varphi_i(X)\cap G^c\subseteq \big(\varphi(X)\cap G^c\big)-\{a\}$. We now remove successively all such elements~$a$ by applying successively all $u_a$ for $a\in  G^c$. It follows that $u_T\circ\varphi$ is a $k$-linear combination of functions $f:X\to T$ such that $f(X)\cap G^c=\emptyset$, that is, $f(X)\subseteq G$.\mpn

(b) Since $\Theta_{T,X}:F_T(X)\to \S_{E,R\op}(X)$ is surjective, $\S_{E,R\op}(X)$ is generated as a $k$-module by the images $\Theta_{T,X}(\varphi)$ of all maps $\varphi:X\to T$. For any $a\in G^c$, $u_a\circ\varphi$ has the same image as $\varphi$ under $\Theta_{T,X}$, by Theorem~\ref{phi minus uaphi}. Therefore $u_T\circ\varphi$ has the same image as $\varphi$ under~$\Theta_{T,X}$.
Moreover, $u_T\circ\varphi$ is a $k$-linear combination of functions $f:X\to T$ such that $f(X)\subseteq G$, by~(a). Finally, if $E\nsubseteq f(X)$, then $f\in \Ker\Theta_{T,X}$ by Theorem~\ref{surjection}, so we can remove any such function in the linear combination $u_T\circ\varphi$ without changing the image $\Theta_{T,X}(u_T\circ\varphi)$.
So we are left with linear combinations of maps $f:X \to T$ such that $E\subseteq f(X) \subseteq G$.
\endpf

We now mention that much more can be said about the elements $u_a$ appearing in Theorem~\ref{generators-SER}.

\result{Definition} \label{graph}
Let $T$ be a finite lattice. Recall that $G^c$ denotes the complement of~$G$ in~$T$. We define an oriented graph structure $\mathcal{G}(T)$ on $T$ in the following way~: for $x,y\in T$, there is an edge $\fleche{x}{y}$ from $x$ to $y$ in $\mathcal{G}(T)$ if there exists $a\in G^c$ such that $(x,y)$ is a pair of consecutive elements in the reduction sequence associated to~$a$.
\fresult

\result{Theorem} \label{idempotents-ua} Keep the notation of Theorem~\ref{generators-SER} and let $\mathcal{G}(T)$ be the graph structure on $T$ introduced in Definition~\ref{graph}.
\begin{enumerate}
\item The graph $\mathcal{G}(T)$ has no (oriented or unoriented) cycles, and each vertex has at most one outgoing edge. Hence $\mathcal{G}(T)$ is a forest.
\item For $a\in G^c$, the element $u_a$ is an idempotent of $k(T^T)$.
\item $u_a\circ u_b=u_b\circ u_a$ for any $a,b\in G^c$.
\item The element $u_T$ is an idempotent of~$k(T^T)$.
\end{enumerate}
\fresult

There is actually a closed formula for~$u_T$ and this is useful for the explicit description of the action of correspondences on the evaluation of simple functors (see Theorem~\ref{action}).
Otherwise, Theorem~\ref{idempotents-ua} has apparently no direct implication on the structure of correspondence functors, so we omit the proof.


\section{Linear independence of the generators} \label{Section-independence}

\medskip
\noindent
In Section~\ref{Section-generators}, we found a set~$\Theta_{T,X}(\CB_X)$ of generators for the evaluation $\S_{E,R\op}(X)$ of a fundamental functor~$\S_{E,R\op}$. We now move to linear independence.

\result{Theorem} \label{independent-SER}
Let $T$ be a finite lattice, let $(E,R)=\Irr(T)$, let $X$ be a finite set, and let $\CB_X$ be the set of all maps $\varphi:X \to T$ such that $E\subseteq \varphi(X) \subseteq G$, where $G=G(T)$ is the subset defined in Notation~\ref{notation-G}.
The elements $\Theta_{T,X}(\varphi)$, for $\varphi\in\CB_X$, are $k$-linearly independent in~$\S_{E,R\op}(X)$.
\fresult

\pf We consider again the kernel of the map
$$\Theta_{T,X}:F_T(X)\longrightarrow \S_{E,R\op}(X) \mvirg$$
which was described in Theorem~\ref{kertheta} by a system of linear equations. This can be reformulated by introducing the $k$-linear map
$$\xymatrix@R=.1ex{
\eta_{E,R,X}:F_T(X)\ar[r]&F_{\Iup(E,R)}(X)\\
{\phantom{\eta_{E,R,X}:}}\varphi\hspace{4ex}\ar@{|->}[r]&*!U(0.7){\sumb{\psi:X\to \Iup(E,R)}{\varphi\vdashER{E,R}\psi}\limits\psi}
}
$$
where the notation $\varphi\vdashER{E,R}\psi$ means, as before, that $\varphi:X\to T$ and $\psi:X\to \Iup(E,R)$ satisfy the equivalent conditions of Theorem~\ref{Gammapsigamma}.
Theorem~\ref{kertheta} asserts that
$$\Ker(\Theta_{T,X})=\Ker(\eta_{E,R,X}) \mpoint$$
For handling the condition $\varphi\vdashER{E,R}\psi$, we shall use part~(e) of Theorem~\ref{Gammapsigamma}, namely
\begin{equation}\label{vdash(e)}
\varphi\vdashER{E,R}\psi\iff 
\left\{\begin{array}{l} \varphi\leq\wedge\psi \mvirg\\
\rule{0ex}{3ex}
\forall e\in E,\;\exists x\in X \;\;\text{such that} \;\; \varphi(x)=e \;\;\text{and} \;\; \psi(x)=[e,\cdot[_E \mpoint
\end{array}\right.
\end{equation}
\par

Let $N=N_{E,R,X}$ be the matrix of $\eta_{E,R,X}$ with respect to the standard basis of $F_T(X)$, consisting of maps $\varphi:X\to T$, and the standard basis of $F_{\Iup(E,R)}(X)$, consisting of maps $\psi:X\to \Iup(E,R)$.
Explicitly,
\begin{equation}\label{matrixN}
N_{\psi,\varphi}=\left\{\begin{array}{ll}1&\hbox{if}\;\varphi\vdashER{E,R}\psi \mvirg\\
0&\hbox{otherwise.}\end{array}\right.
\end{equation}
Note that $N$ is a square matrix in the special case when $T=\Idown(E,R)$, because complementation yields a bijection between $\Idown(E,R)$ and~$\Iup(E,R)$. However, if $T$ is a proper quotient of~$\Idown(E,R)$, then $N$ has less columns.\par

In order to prove that the elements $\Theta_{T,X}(\varphi)$, for $\varphi\in\CB_X$, are $k$-linearly independent, we shall prove that the elements $\eta_{E,R,X}(\varphi)$, for $\varphi\in\CB_X$, are $k$-linearly independent.
In other words, we have to show that the columns of~$N$ indexed by $\varphi\in\CB_X$ are $k$-linearly independent. Now we consider only the rows indexed by elements of the form $\psi=\zeta\circ\varphi'$, where $\varphi'\in\CB_X$ and $\zeta:G\to \Iup(E,R)$ is the map defined in (\ref{def-zeta}).
We then define the square matrix~$M$, indexed by $\CB_X\times\CB_X$, by
$$\forall \varphi,\varphi'\in\CB_X,\;\;M_{\varphi',\varphi}=N_{\zeta\circ\varphi',\varphi} \mpoint$$
We are going to prove that $M$ is invertible and this will prove the required linear independence.\par

The invertibility of~$M$ implies in particular that the map~$\zeta$ must be injective, otherwise two rows of~$M$ would be equal. Therefore $M$ turns out to be a submatrix of~$N$, but this cannot be seen directly from its definition (unless an independent proof of the injectivity of~$\zeta$ is provided).\par

The characterization of the condition $\varphi\vdashER{E,R}\psi$ given in~(\ref{vdash(e)}) implies that
$$M_{\varphi',\varphi}=\left\{\begin{array}{ll}
1&\hbox{if}\;\varphi\leq\wedge\zeta\varphi'\; \hbox{and}\;
\forall e\in E,\exists x\in X,\;\varphi(x)=e=\varphi'(x) \mvirg\\
0&\hbox{otherwise}\mvirg\end{array}\right.$$
because the equality $\zeta\varphi'(x)=[e,\cdot[_E$ is equivalent to $\varphi'(x)=e$, by definition of~$\zeta$
(see Notation~\ref{def-zeta}).\par

By Proposition~\ref{rho-leq}, if $t,t'\in G$ are such that $t\leq\wedge \zeta(t')$, then $r^\infty(t)\leq r^\infty(t')$ and $\sigma^\infty(t)\leq \sigma^\infty(t')$. Let $\preceq$ be the preorder on $G$ defined by these conditions, i.e. for all $t,t'\in G$,
$$t\preceq t' \iff r^\infty(t)\leq r^\infty(t')\;\hbox{and}\;\sigma^\infty(t)\leq \sigma^\infty(t')\mpoint$$
We extend this preorder to $\CB_X$ by setting, for all $\varphi', \varphi\in\CB_X$,
$$\varphi\preceq\varphi' \iff \;\forall x\in X,\;\varphi(x)\preceq \varphi'(x)\mvirg$$
which makes sense because $\varphi(x),\varphi'(x)\in G$ by definition of~$\CB_X$.
We denote by $\equi$ the equivalence relation defined by this preorder.\par

Clearly the condition $M_{\varphi',\varphi}\neq 0$ implies $\varphi\leq\wedge\zeta\varphi'$, hence $\varphi\preceq\varphi'$ by Proposition~\ref{rho-leq} quoted above.
In other words the matrix $M$ is block triangular, the blocks being indexed by the equivalence  classes of the preorder $\preceq$ on $\CB_X$. Showing that $M$ is invertible is equivalent to showing that all its diagonal blocks are invertible. In other words, we must prove that, for each equivalence class $C$ of~$\CB_X$ for the relation~$\equi$, the matrix $M_C=(M_{\varphi',\varphi})_{\varphi',\varphi\in C}$ is invertible. Let $C$ be such a fixed equivalence class.\mpn

Recall from Definition~\ref{notation-G} that $G=\Lambda E\sqcup \widehat G$.
If $t\in \widehat G$, then by Lemma~\ref{bulb} $e:=\sigma^\infty(t)$ belongs to~$E$ and $[t,e]= \{r^k(e),r^{k-1}(e), \ldots, r^1(e),e \}$, where $t=r^k(e)=r^\infty(e)$.
By Lemma~\ref{interval}, all elements of $[\,t,\sigma^\infty(t)\,]_T$ belong to~$E$ except $t$ itself.
Moreover, if $x\in [\,t,\sigma^\infty(t)\,]_T$ then $r^\infty(x)=t$ by Lemma~\ref{bulb}.
It follows that the sets
$$G_t=[\,t,\sigma^\infty(t)\,]_T \,, \quad \text{for } t\in \widehat G \mvirg$$
are disjoint, and contained in~$G$. Let
$$G_*=G-\bigsqcup_{t\in \widehat G}G_t\mvirg$$
so that we get a partition
$$G=\bigsqcup_{t\in\{*\}\sqcup \widehat G} G_t \mpoint$$

\result{Lemma} \label{partition} Let $\varphi',\varphi\in\CB_X$. If $\varphi'\equi\varphi$, then for all $t\in \{*\}\sqcup \widehat G$,
$$\varphi'^{-1}(G_t)=\varphi^{-1}(G_t)\mpoint$$
\fresult

\pf Let $t\in \widehat G$ and $x\in \varphi^{-1}(G_t)$. Then $\varphi(x)\in[\,t,\sigma^\infty(t)\,]_T$,
hence $r^\infty\varphi(x)=t$ and $\sigma^\infty\varphi(x)=\sigma^\infty(t)$, by Lemma~\ref{bulb}.
But the relation $\varphi'\equi\varphi$ implies that $r^\infty\varphi'(x)=r^\infty\varphi(x)$ and $\sigma^\infty\varphi'(x)=\sigma^\infty\varphi(x)$. Therefore $r^\infty\varphi'(x)=t$ and $\sigma^\infty\varphi'(x)=\sigma^\infty(t)$, from which it follows that $\varphi'(x)\in[\,t,\sigma^\infty(t)\,]_T$, that is, $x\in \varphi'^{-1}(G_t)$. This shows that $\varphi^{-1}(G_t)\subseteq\varphi'^{-1}(G_t)$. By exchanging the roles of $\varphi$ and $\varphi'$, we obtain $\varphi'^{-1}(G_t)=\varphi^{-1}(G_t)$.\par

Now $G_*$ is the complement of $\bigsqcup\limits_{t\in \widehat G}G_t$ in~$G$ and the functions $\varphi',\varphi$ have their values in~$G$ (by definition of~$\CB_X$).
So we must have also $\varphi'^{-1}(G_*)=\varphi^{-1}(G_*)$.
\endpf

For every $t\in\{*\}\sqcup \widehat G$, we define
$$X_t=\varphi_0^{-1}(G_t)$$
where $\varphi_0$ is an arbitrary element of~$C$.
It follows from Lemma~\ref{partition} that this definition does not depend on the choice of~$\varphi_0$.
Therefore, the equivalence class~$C$ yields a partition
$$X=\bigsqcup_{t\in\{*\}\sqcup \widehat G} X_t \mvirg$$
and every function $\varphi\in C$ decomposes as the disjoint union of the functions $\varphi_t$,
where $\varphi_t:X_t\to G_t$ is the restriction of $\varphi$ to~$X_t$.\par

For $t\in \widehat G$, define
$$E_t=]\,t,\sigma^\infty(t)\,]_T \mpoint$$
By Lemma~\ref{interval}, this consists of elements of~$E$, so $E_t=E\cap G_t$.
Then we define $E_*=E-\bigsqcup\limits_{t\in \widehat G}E_t=E\cap G_*$, so that we get a partition
$$E=\bigsqcup_{t\in\{*\}\sqcup \widehat G} E_t \mpoint$$
For every $t\in\{*\}\sqcup \widehat G$ and for every $\varphi\in C$, the fact that $\varphi$ belongs to~$\CB_X$ implies that $\varphi_t$ belongs to the set~$\CB_{X,t}$ of all maps $\varphi_t:X_t\to G_t$ such that
$$E_t\subseteq \varphi_t(X_t)\subseteq G_t \mpoint$$
Moreover, if $\varphi', \varphi\in C$, then
$$M_{\varphi',\varphi}=1\iff \forall t\in\{*\}\sqcup \widehat G,\;\left\{\begin{array}{l}\forall x \in X_t,\;\varphi_t(x)\leq\meet\zeta\varphi'_t(x)\\\rule{0ex}{3ex} \forall e\in E_t,\;\exists x\in X_t,\;\varphi'(x)=e=\varphi(x)\mpoint\end{array}\right.$$
It follows that the matrix $M_C$ is the tensor product of the square matrices $M_{C,t}$ for $t\in \{*\}\sqcup \widehat G$, where
the matrix $M_{C,t}$ is indexed by the functions $\varphi_t:X_t\to G_t$ in~$\CB_{X,t}$ and satisfies
$$(M_{C,t})_{\varphi'_t,\varphi_t}=1\iff \left\{\begin{array}{l}\forall x \in X_t,\;\varphi_t(x)\leq\meet\zeta\varphi'_t(x)\\\rule{0ex}{3ex} \forall e\in E_t,\;\exists x\in X_t,\;\varphi'(x)=e=\varphi(x)\mpoint\end{array}\right.$$
In order to show that $M_C$ is invertible, we shall prove that each matrix $M_{C,t}$ is invertible.\mpn

If $\varphi_*\in\CB_{X,*}$ and $x\in X_*$, then $\varphi_*(x)\in G_*$, hence $\varphi_*(x)\notin \widehat G$, because $\widehat G$ consists of all the bottom elements of the intervals~$G_t$ where $t\in \widehat G$. Therefore, the condition $\varphi_*(x)\leq \meet\zeta\varphi'_*(x)$ implies $\varphi_*(x)\leq\varphi'_*(x)$ by Proposition~\ref{rho-leq}, because $\varphi'_*(x)\notin \widehat G$. It follows that the matrix $M_{C,*}$ is unitriangular, hence invertible, as required.\par

Now we fix $t\in \widehat G$, we consider the matrix $M_{C,t}$ and we discuss the special role played by the elements of the set~$\widehat G$.
The interval $G_t=[t,\sigma^{\infty}(t)]_T$ is isomorphic to the totally ordered lattice~$\sou{n}=\{0,1,\ldots,n\}$, for some~$n\geq1$, and the set of irreducible elements $E_t=]t,\sigma^{\infty}(t)]_T$ is isomorphic to~$[n]=\{1,\ldots,n\}$.
Composing the maps $\varphi_t:X_t\to G_t$ with this isomorphism, we obtain maps $X_t\to \sou{n}$.\par

Changing notation for simplicity, we write $X$ for $X_t$ and $\varphi$ for~$\varphi_t$, and we let $\CB_X^{\sou{n}}$ be the set of all maps $\varphi:X\to \sou{n}$ corresponding to maps in~$\CB_{X,t}$, i.e. satisfying the condition $[n]\subseteq \varphi(X)\subseteq\sou{n}$.
We note that this condition is the same as the condition that $\varphi$ belongs to~$\CB_X$ for the lattice~$\sou{n}$, because the set $G(\sou{n})$ is the whole of~$\sou{n}$ by Example~\ref{Example-tot}.
The matrix $M_{C,t}$, which we write~$M^{\sou{n}}$ for simplicity, is now indexed by all the maps in~$\CB_X^{\sou{n}}$ and we have
$$M_{\varphi',\varphi}^{\sou{n}}=1\iff \left\{\begin{array}{l}\forall x \in X,\;\varphi(x)\leq\meet\zeta\varphi'(x)\\
\rule{0ex}{3ex} \forall e\in [n],\;\exists x\in X,\;\varphi'(x)=e=\varphi(x)\mpoint\end{array}\right.$$
Here we need to clarify the meaning of the notation $\meet\zeta$, so we recall that for any $g\in G$, we have
$$\wedge\zeta(g)=\left\{\begin{array}{ll}g&\hbox{if}\; g\in E \mvirg\\
\sigma^\infty(g)&\hbox{if}\;g\notin E\mpoint \end{array}\right.$$
If $g$ belongs to~$G_t$ and is mapped to~$a\in\sou{n}$ via the isomorphism $G_t\cong\sou{n}$, then $\sigma^\infty(g)$ is mapped to~$n$ and we obtain
$$\wedge\zeta(a)=\left\{\begin{array}{ll}a&\hbox{if}\;a\in[n], \; \text{i.e. } a\neq 0 \mvirg\\
n&\hbox{if}\;a=0\mpoint \end{array}\right.$$
The point here is that we obtain the same result as the one we would have obtained by working with the lattice~$\sou{n}$, that is, by working with the corresponding map $\zeta:\sou{n}\to \Iup([n],\tot)$, which is easily seen to be a bijection, mapping $0$ to $\emptyset$ and $j\geq1$ to~$[j,n]$.\par

Now we return to the beginning of the proof of Theorem~\ref{independent-SER} in the special case of the lattice $\sou{n}$, with $\Irr(\sou{n})=([n],\tot)$, where $\tot$ denotes the usual total order.
We have
$$\Ker(\Theta_{\sou{n},X})=\Ker(\eta_{[n],\tot,X})$$
and the matrix $N$ of~$\eta_{[n],\tot,X}$ has entries 0 and~1, with
$$N_{\psi,\varphi}=1\iff \left\{\begin{array}{l}\forall x \in X,\;\varphi(x)\leq\meet\psi(x)\\
\rule{0ex}{3ex} \forall e\in [n],\;\exists x\in X,\;\varphi(x)=e \;\text{ and } \;\psi(x)=[e,\cdot[_{[n]}\mvirg\end{array}\right.$$
where $\varphi:X\to\sou{n}$ and $\psi:X\to \Iup([n],\tot)$.
But since $\zeta:\sou{n}\to \Iup([n],\tot)$ is a bijection, we can write $\psi=\zeta\varphi'$ and index the rows by the set of all functions $\varphi':X\to\sou{n}$.
We obtain
$$N_{\varphi',\varphi}=1\iff \left\{\begin{array}{l}\forall x \in X,\;\varphi(x)\leq\meet\zeta\varphi'(x)\\
\rule{0ex}{3ex} \forall e\in [n],\;\exists x\in X,\;\varphi'(x)=e=\varphi(x)\mvirg\end{array}\right.$$
If $\varphi$ or $\varphi'$ is not in~$\CB_X^{\sou{n}}$ (that is, the image of either $\varphi$ or $\varphi'$ does not contain $[n]$), then the second condition cannot hold and so $N_{\varphi',\varphi}=0$. Thus we restrict the matrix~$N$ to the rows and columns indexed by~$\CB_X^{\sou{n}}$.
This restriction is exactly the same matrix as the matrix~$M^{\sou{n}}$ above.
Therefore, in order to prove that $M^{\sou{n}}$ is invertible, it suffices to prove that the columns of~$N$ indexed by~$\CB_X^{\sou{n}}$ are $k$-linearly independent.
This in turn is equivalent to the condition that the set 
$\eta_{[n],\tot,X}(\CB_X^{\sou{n}})$ is $k$-linearly independent, or also that the set
$\Theta_{\sou{n},X}(\CB_X^{\sou{n}})$ is $k$-linearly independent in~$\S_{[n],\tot\op}(X)$, because $\Ker(\Theta_{\sou{n},X})=\Ker(\eta_{[n],\tot,X})$.
In other words, we have to prove Theorem~\ref{independent-SER} in the case of a total order.\par

By Theorem~\ref{surjection}, we know that the surjective morphism
$$\Theta_{\sou{n}}: F_{\sou{n}} \longrightarrow \S_{[n],\tot\op}$$
has $H_{\sou{n}}$ in its kernel, where $H_{\sou{n}}$ is the subfunctor of~$F_{\sou{n}}$ generated by all the maps $\varphi:X\to \sou{n}$ such that $[n]\not\subseteq \varphi(X)$.
Therefore $\Theta_{\sou{n}}$ induces a surjective morphism
$$\overline\Theta_{\sou{n}}: F_{\sou{n}}/H_{\sou{n}} \longrightarrow \S_{[n],\tot\op} \mpoint$$
But Theorem 11.8 in~\cite{BT3} asserts that $F_{\sou{n}}/H_{\sou{n}}$ is isomorphic to~$\S_{[n],\tot}$, hence also to $\S_{[n],\tot\op}$ in view of the poset isomorphism $([n],\tot\op)\cong([n],\tot)$ via the map $j\mapsto n-j+1$.
Clearly the set $\CB_X^{\sou{n}}$ is a $k$-basis of $F_{\sou{n}}(X)/H_{\sou{n}}(X)$, so that $\S_{[n],\tot\op}(X)$ is also a free $k$-module of rank~$|\CB_X^{\sou{n}}|$.
Evaluation at~$X$ yields a surjective homomorphism
$$\overline\Theta_{\sou{n},X} : F_{\sou{n}}(X)/H_{\sou{n}}(X) \longrightarrow \S_{[n],\tot\op}(X)$$
between two free $k$-modules of the same rank, hence an isomorphism (by standard algebraic $K$-theory, see Lemma~6.8 in~\cite{BT3}).
Since the elements $\Theta_{\sou{n},X}(\varphi)$, for $\varphi\in\CB_X^{\sou{n}}$, are the images under $\overline\Theta_{\sou{n},X}$ of the $k$-basis $\CB_X^{\sou{n}}$ of~$F_{\sou{n}}(X)/H_{\sou{n}}(X)$, they form a $k$-basis of~$\S_{[n],\tot\op}(X)$. In particular, they are $k$-linearly independent.\par

This completes the proof of Theorem~\ref{independent-SER}.
\endpf

In order to obtain formulas for the dimension of the evaluation of a fundamental functor, we need a combinatorial lemma, which is Lemma~8.1 in~\cite{BT2}.

\result{Lemma} \label{sandwich} Let $E$ be a subset of a finite set~$G$.
For any finite set $X$, the number $N$ of all maps ${\varphi:X\to G}$ such that $E\subseteq \varphi(X)\subseteq G$ is equal to
$$N=\sum_{i=0}^{|E|}(-1)^i{|E|\choose i}(|G|-i)^{|X|}\mpoint$$
\fresult

We can now prove one of our main results about fundamental correspondence functors.
This generalizes a formula obtained in~\cite{BT3} in the case of a total order.

\result{Theorem} \label{dim-SER}
Let $(E,R)$ be a finite poset and let $T$ be any lattice such that $(E,R)=\Irr(T)$.
Let $X$ be a finite set and let $\CB_X$ be the set of all maps $\varphi:X \to T$ such that $E\subseteq \varphi(X) \subseteq G$, where $G=G(T)$ is the subset defined in Notation~\ref{notation-G}.
\begin{enumerate}
\item The set $\Theta_{T,X}(\CB_X)$ (more precisely, the injective image of~$\CB_X$ under~$\Theta_{T,X}$) is a $k$-basis of~$\S_{E,R\op}(X)$.
\item The $k$-module $\S_{E,R\op}(X)$ is free of rank
$${\rm rk}_k\big(\S_{E,R\op}(X)\big)=|\CB_X|=\sum_{i=0}^{|E|}(-1)^i\binom{|E|}{i}(|G|-i)^{|X|}\mpoint$$
\end{enumerate}
\fresult

\pf (a) follows from Theorem~\ref{generators-SER} and Theorem~\ref{independent-SER}.\mpn

(b) The formula follows immediately from (a) and Lemma~\ref{sandwich}.
\endpf

\result{Corollary} \label{choice-G}
With the notation above, $|G|$ only depends on~$(E,R)$, and not on the choice of~$T$.
\fresult

\pf
The formula of Theorem~\ref{dim-SER} implies that
$${\rm rk}_k\big(\S_{E,R\op}(X)\big)\sim |G|^{|X|} \qquad\text{ as } |X|\to\infty\mpoint$$
Since $\S_{E,R\op}$ only depends on~$(E,R)$, it follows that $|G|$ only depends on~$(E,R)$.
\endpf

We shall prove in a future paper a stronger property~: the full subposet $G$ of~$T$ only depends on~$(E,R)$ up to isomorphism.


\section{From fundamental functors to simple functors} \label{Section-simple-tensor}

\medskip
\noindent
In this section, we complete the description of simple functors by showing that they can be constructed directly from fundamental functors.
This uses in an essential way the fact, proved in Theorem~\ref{SER-free} below, that each evaluation of a fundamental functor $\S_{E,R}$ is a free $k\Aut(E,R)$-module, thanks to our main Theorem~\ref{dim-SER}.\par

We first recall some basic constructions of functors (see~\cite{BT2, BT3}).
Let $E$ be a finite set and let $W$ be a left $k\CC(E,E)$-module.
The correspondence functor $L_{E,W}$ is defined by
$$L_{E,W}(X):=k\CC(X,E)\otimes_{k\CC(E,E)}W \mvirg$$
with an obvious left action of correspondences in $\CC(Y,X)$ by composition.
There is a subfunctor $J_{E,W}$ of~$L_{E,W}$ defined as follows (see Lemma~2.5 in~\cite{BT2})~:
$$J_{E,W}(X)=\Big\{\sum_i\varphi_i\otimes w_i \in L_{E,W}(X) \,\mid\, \forall \rho\in k\CC(E,X) \,, \sum_i (\rho\varphi_i)\cdot w_i=0 \Big\} \mpoint$$
Let us mention an important property of the quotient functor $L_{E,W}/J_{E,W}$.

\result{Lemma} \label{Ker-Res}
Let $E$ be a finite set and let $W$ be a left $k\CC(E,E)$-module.
\begin{enumerate}
\item $J_{E,W}(E) =\{0\}$ and $L_{E,W}(E)/J_{E,W}(E)\cong L_{E,W}(E) \cong W$.
\item Let $\alpha \in L_{E,W}(X)/J_{E,W}(X)$ where $X$ is some finite set.
Then $\rho\cdot\alpha=0$ for every $\rho\in k\CC(E,X)$ if and only if $\alpha=0$.
\end{enumerate}
\fresult

\pf
(a) It is clear that $L_{E,W}(E) = k\CC(E,E)\otimes_{k\CC(E,E)}W\cong W$. Corresponding to $w\in W$, let $\Id\otimes w\in L_{E,W}(E)$.
If $\Id\otimes w\in J_{E,W}(E)$, we choose $\rho=\Id\in \CC(E,E)$ and we obtain $w=(\rho\circ\Id)\cdot w=0$, by definition of~$J_{E,W}(X)$.
This shows that $J_{E,W}(E) =\{0\}$. \mpn

(b) Let $\rho\in k\CC(E,X)$. It follows from (a) that there is a commutative diagram
$$\xymatrix{
0 \ar[r] & J_{E,W}(X) \ar[r] \ar[d] & L_{E,W}(X) \ar[r]^-{\pi} \ar[d]^-{\rho} & L_{E,W}(X)/J_{E,W}(X) \ar[r] \ar[d]^-{\rho} & 0 \\
& \;0\; \ar[r] & L_{E,W}(E) \ar[r]^-{\cong} & W \ar[r] & 0
}$$
Let $\sum_i\varphi_i\otimes w_i \in L_{E,W}(X)$ such that $\pi(\sum_i\varphi_i\otimes w_i )=\alpha$.
From the assumption that $\rho\cdot\alpha=0$ for every $\rho\in k\CC(E,X)$, we obtain
$$0=\rho\cdot \big(\sum_i\varphi_i\otimes w_i \big)=\sum_i\rho\varphi_i\otimes w_i \in L_{E,W}(E) \mpoint$$
Viewing this in~$W$ via the isomorphism $L_{E,W}(E) \cong W$, we get $\sum_i (\rho\varphi_i)\cdot w_i=0$ for every $\rho\in k\CC(E,X)$.
In other words, $\sum_i\varphi_i\otimes w_i \in J_{E,W}(X)$ and it follows that $\alpha=0$.
\endpf

We now recall the construction of fundamental functors and simple functors (see~\cite{BT2, BT3}), which are special cases of the construction above.
Using the fundamental $k\CC(E,E)$-module $\CP_E f_R$ associated with a poset $(E,R)$, we obtain the fundamental functor
$$\S_{E,R} := L_{E,\CP_E f_R} / J_{E,\CP_E f_R} \mpoint$$
If now $V$ is a left $k\Aut(E,R)$-module, we define the $k\CC(E,E)$-module
$$T_{R,V} :=\CP_E f_R\otimes_{k\Aut(E,R)}V$$
(using the right $k\Aut(E,R)$-module structure on~$\CP_E f_R$ described in Proposition~\ref{fundamental-module}) and we obtain an associated correspondence functor
$$S_{E,R,V} := L_{E,T_{R,V}} / J_{E,T_{R,V}} \mpoint$$
When $k$ is a field and $V$ is a simple $k\Aut(E,R)$-module, this provides the explicit construction of the simple correspondence functor $S_{E,R,V}$ (appearing in Proposition~\ref{SERV} and Theorem~\ref{parametrization}).\par

Using the right $k\Aut(E,R)$-module structure on~$\CP_E f_R$, we can define a right $k\Aut(E,R)$-module structure on each evaluation
$$L_{E,\CP_E f_R}(X)=k\CC(X,E)\otimes_{k\CC(E,E)} \CP_E f_R$$
and we now show that this right module structure can be carried to~$\S_{E,R}(X)$.

\result{Lemma} \label{Aut-action}
\begin{enumerate}
\item $J_{E,\CP_E f_R}(X)$ is a right $k\Aut(E,R)$-submodule of~$L_{E,\CP_E f_R}(X)$.
\item $\S_{E,R}(X)$ has a right $k\Aut(E,R)$-module structure.
\item The left action of any element of $k\CC(Y,X)$ is a homomorphism of right $k\Aut(E,R)$-modules $\S_{E,R}(X)\to \S_{E,R}(Y)$.
\end{enumerate}
\fresult

\pf
(a) Since $\CP_E$ is a quotient algebra of $k\CC(E,E)$ and the tensor product defining $L_{E,\CP_E f_R}(X)$ is over $k\CC(E,E)$,
any element of~$L_{E,\CP_E f_R}(X)$ can be written $\varphi\otimes f_R$ for some $\varphi\in k\CC(X,E)$.
By Proposition~\ref{fundamental-module}, the right action of $\sigma\in\Aut(E,R)$ is given by
$$(\varphi\otimes f_R)\Delta_\sigma =\varphi\otimes (\Delta_\sigma f_R) =(\varphi \Delta_\sigma)\otimes f_R \mpoint$$
If $\varphi\otimes f_R\in J_{E,\CP_E f_R}(X)$, then $(\rho\varphi)\cdot f_R=0$ for all $\rho\in k\CC(E,X)$.
Then the element $(\varphi \Delta_\sigma)\otimes f_R$ satisfies
$$(\rho\varphi \Delta_\sigma)\cdot f_R=(\rho\varphi )\cdot f_R\Delta_\sigma=0$$
for all $\rho\in k\CC(E,X)$. Therefore $(\varphi \Delta_\sigma)\otimes f_R=(\varphi \otimes f_R)\Delta_\sigma$ belongs to~$J_{E,\CP_E f_R}(X)$,
as was to be shown.\mpn

(b) This follows immediately from~(a).\mpn

(c) This follows from the fact that the left and right actions commute,
by associativity of the composition $k\CC(Y,X)\times k\CC(X,E) \times k\CC(E,E) \to k\CC(Y,E)$.
\endpf

Given a finite lattice~$T$ with $\Irr(T)=(E,R)$, we can now explain how the morphism $\Theta_T: F_T \to \S_{E,R\op}$ of Theorem~\ref{surjection} is defined.
This appears in Theorem~6.5 of~\cite{BT3}.

\result{Definition} \label{Def-Theta}
Let $\varphi:X\to T$ be a map, i.e. a generator of~$F_T(X)$. Then $\Theta_{T,X}(\varphi)$ is the class in the quotient~$\S_{E,R\op}(X)$ of the element
$$\Gamma_\varphi\otimes f_{R\op} \in k\CC(X,E)\otimes_{k\CC(E,E)} \CP_E f_{R\op} = L_{E,\CP_E f_{R\op}}(X) \mvirg$$
where $\Gamma_\varphi$ is the correspondence defined in Notation~\ref{Gamma}.
\fresult

In order to use the action of automorphisms, we first need a lemma.

\result{Lemma} \label{Aut(T)}
Let $(E,R)$ be a finite poset.
\begin{enumerate}
\item For any finite lattice $T$ such that $\Irr(T)= (E,R)$, restriction induces an injective group homomorphism $\Aut(T)\to\Aut(E,R)$.
\item There exists a finite lattice $T$ such that $\Irr(T)= (E,R)$ and such that the restriction homomorphism $\Aut(T)\to\Aut(E,R)$ is an isomorphism.
\end{enumerate}
\fresult

\pf
(a) For any lattice $T$ such that $\Irr(T)= (E,R)$, any lattice automorphism of $T$ induces an automorphism of the poset $(E,R)$. This gives a group homomorphism $\Aut(T)\to\Aut(E,R)$, which is injective since any element $t$ of $T$ is equal to the join of the irreducible elements $e\leq_Tt$.\mpn

(b) Requiring that $\Aut(T)\cong\Aut(E,R)$ amounts to requiring that any automorphism of $(E,R)$ can be extended to an automorphism of~$T$. This is clearly possible if we choose $T=\Idown(E,R)$.
\endpf

The group $\Aut(T)$ acts on the right on~$F_T(X)$ as follows~:
$$\forall \varphi: X\to T \,, \; \forall \sigma\in\Aut(T) \,, \qquad \varphi\cdot \sigma := \sigma^{-1}\circ\varphi \mpoint$$

\result{Lemma} \label{equivariant-Theta} Let $(E,R)$ be a poset.
Let $T$ be a lattice such that $\Irr(T)= (E,R)$ and such that the restriction homomorphism $\Aut(T)\to\Aut(E,R)$ is an isomorphism.
For every finite set~$X$, the map
$$\Theta_{T,X}: F_T(X) \longrightarrow \S_{E,R\op}(X)$$
is a homomorphism of right $k\Aut(E,R)$-modules.
\fresult

\pf
First note that we obviously have an equality $\Aut(E,R)=\Aut(E,R\op)$.
Let us first prove that, for any $\varphi:X\to T$ and any $\sigma\in \Aut(T)\cong \Aut(E,R)$,
$$\Gamma_\varphi\Delta_\sigma=\Gamma_{\sigma^{-1}\varphi} \mpoint$$
An element $(x,e)\in X\times E$ belongs to the left hand side if and only if $(x,\sigma(e))\in\Gamma_\varphi$, because $(\sigma(e),e)\in \Delta_\sigma$. 
By definition of~$\Gamma_\varphi$, this is equivalent to the condition $\sigma(e)\leq_T \varphi(x)$, which in turn is equivalent to $e\leq_T \sigma^{-1}\varphi(x)$ because $\sigma\in\Aut(T)$. Thus we obtain that $(x,e)$ satisfies the condition defining~$\Gamma_{\sigma^{-1}\varphi}$, that is, $(x,e)$ belongs to the right hand side.\par

Now we can compute $\Theta_{T,X}(\varphi\cdot\sigma)=\Theta_{T,X}(\sigma^{-1}\varphi)$.
By Definition~\ref{Def-Theta}, this is the class in the quotient~$\S_{E,R\op}(X)$ of the element
$$\Gamma_{\sigma^{-1}\varphi} \otimes f_{R\op} = \Gamma_\varphi\Delta_\sigma\otimes f_{R\op}
=\Gamma_\varphi\otimes \Delta_\sigma f_{R\op}=(\Gamma_\varphi\otimes f_{R\op})\Delta_\sigma \mvirg$$
using the definition of the right action of~$\sigma$ on~$\Gamma_\varphi\otimes f_{R\op}$ (Proposition~\ref{fundamental-module}). 
Since the class of~$\Gamma_\varphi\otimes f_{R\op}$ is $\Theta_{T,X}(\varphi)$,
this shows that $\Theta_{T,X}(\varphi\cdot\sigma)=\Theta_{T,X}(\varphi)\cdot\sigma$, as required. \mpn
\endpf

Our next result is the key for the rest of this section, hence for the computation of the dimension of the evaluations of simple functors (Theorem~\ref{dim-SERV}).

\result{Theorem} \label{SER-free} Let $(E,R)$ be a poset.
Let $T$ be a lattice such that $\Irr(T)= (E,R)$ and such that the restriction homomorphism $\Aut(T)\to\Aut(E,R)$ is an isomorphism.
For every finite set~$X$, the evaluation $\S_{E,R\op}(X)$ is a free right $k\Aut(E,R)$-module.
\fresult

\pf
The set of all maps $\varphi:X\to T$ is a $k$-basis of~$F_T(X)$ and is permuted by the right action of~$\Aut(E,R)$.
If $G$ is as before (see Theorem~\ref{dim-SER}), 
we claim that the subset $\CB_X$ of all maps satisfying $E\subseteq \varphi(X) \subseteq G$ is freely permuted by~$\Aut(E,R)$.
First note that $\Aut(E,R)$ obviously leaves $E$ invariant. It also leaves $G$ invariant because $\Aut(T)\cong \Aut(E,R)$ preserves the characterization of~$G$ given in Lemma~\ref{characterize-G}. Therefore $\Aut(E,R)$ acts on~$\CB_X$.\par

If $\sigma\in\Aut(E,R)$ stabilizes some $\varphi\in\CB_X$, that is, $\varphi\cdot\sigma=\varphi$, then $\sigma^{-1}\varphi(x)=\varphi(x)$ for all $x\in X$, hence in particular $\sigma^{-1}(e)=e$ for every $e\in E$ because $E\subseteq \varphi(X)$. It follows that $\sigma$ is the identity automorphism of~$E$.
This proves the claim above.\par

Now $\CB_X$ is mapped bijectively onto $\Theta_{T,X}(\CB_X)$, which is a $k$-basis of~$\S_{E,R\op}(X)$, by Theorem~\ref{dim-SER}.
By Lemma~\ref{equivariant-Theta}, $\Theta_{T,X}$ is a homomorphism of $k\Aut(E,R)$-modules.
It follows that $\Theta_{T,X}(\CB_X)$ is freely permuted by~$\Aut(E,R)$ and therefore $\S_{E,R\op}(X)$ is a free (right) $k\Aut(E,R)$-module.
\endpf

Now we construct a morphism $\Psi: \S_{E,R}\otimes_{k\Aut(E,R)}V \to S_{E,R,V}$, which will be proved later to be an isomorphism (Theorem~\ref{iso-SERV}).

\result{Proposition} \label{diagram}
Let $(E,R)$ be a finite poset, let $A=\Aut(E,R)$, and let $V$ be a left $kA$-module, generated by a single element~$v$ (e.g. a simple module). 
\begin{enumerate}
\item For any finite set~$X$, there is a commutative diagram
$$\xymatrix{
0 \ar[r] & J_{E,\CP_E f_R}(X) \ar[r]^-{j} \ar[d] & L_{E,\CP_E f_R}(X) \ar[r]^-{\pi} \ar[d] & \S_{E,R}(X) \ar[r] \ar[d]^-{\Id\otimes v} & 0 \\
& J_{E,\CP_E f_R}(X)\otimes_{kA} V \ar[r]^-{j\otimes\Id_V} \ar[d] &
L_{E,\CP_E f_R}(X)\otimes_{kA} V \ar[r]^-{\pi\otimes\Id_V} \ar[d]^-{\Id} & \S_{E,R}(X)\otimes_{kA}V \ar[r] \ar[d]^-{\Psi_X} & 0 \\
0 \ar[r] & J_{E,\CP_E f_R\otimes V}(X) \ar[r]^-{i} & L_{E,\CP_E f_R\otimes V}(X) \ar[r] & S_{E,R,V}(X) \ar[r] & 0
}
$$
\item On the right hand side, both maps $\Id\otimes v$ and $\Psi_X$ are surjective.
\end{enumerate}
\fresult

\pf
The first row comes from the definition of~$\S_{E,R}$ and $j$ denotes the inclusion map while $\pi$ is the quotient map.
The second row is obtained from the first by tensoring with~$V$ (tensoring is right exact), using the right $kA$-module structure obtained in Lemma~\ref{Aut-action}.
The three vertical maps from the first to the second row are all given by $\alpha\mapsto \alpha\otimes v$ and they are surjective because $v$ generates~$V$, hence $V=kA\cdot v$.
The third row comes from the definition of~$S_{E,R,V}$ and $i$ denotes the inclusion map.
Now we have to describe the vertical maps from the second to the third row.
The middle vertical map is the identity because
$$L_{E,\CP_E f_R}(X)\otimes_{kA} V= k\CC(X,E)\otimes_{k\CC(E,E)}\CP_E f_R\otimes_{kA} V
=L_{E,\CP_E f_R\otimes_{kA} V}(X) \mpoint$$
We claim that $(j\otimes\Id_V)\big(J_{E,\CP_E f_R}(X)\otimes_{kA} V\big)$
is contained in~$J_{E,\CP_E f_R\otimes V}(X)$.
It will then follow that $j\otimes\Id_V$ defines the vertical map on the left.
This in turn shows that $\Id$ induces the vertical map~$\Psi_X$ on the right and $\Psi_X$ is surjective.\par

In order to prove the claim, let $\varphi\otimes f_R\in J_{E,\CP_E f_R}(X)$.
This means that
$$\forall \rho\in k\CC(E,X) \,, \qquad(\rho\varphi)\cdot f_R=0 \mpoint$$
It follows that $(\rho\varphi)\cdot (f_R\otimes v)=0$ in~$\CP_E f_R\otimes_{kA} V$ because $\rho\varphi$ only acts on the first term of the tensor product.
This means that the element
$$\varphi\otimes (f_R\otimes v) \in L_{E,\CP_E f_R\otimes V}(X)$$
actually belongs to~$J_{E,\CP_E f_R\otimes V}(X)$. But this element is equal to
$$(\varphi\otimes f_R)\otimes v=(j\otimes\Id_V)(\varphi\otimes f_R\otimes v) \mvirg$$
proving the claim.
\endpf

\result{Notation} \label{Notation-Psi}
Consider the diagram of Proposition~\ref{diagram}.
When $X$ is allowed to vary, the morphisms $\Psi_X$ on the right hand side define a surjective morphism of correspondence functors
$$\Psi:\S_{E,R}\otimes_{kA} V \longrightarrow S_{E,R,V} \mvirg$$
providing a direct link between the fundamental functor $\S_{E,R}$ and the simple functor $S_{E,R,V}$ when $V$ is simple.\par

Similarly, the right hand side composition $\Psi_X\circ(\Id\otimes v)$ yields a surjective morphism
$$\Phi :=\Psi\circ(\Id\otimes v): \S_{E,R} \longrightarrow S_{E,R,V}$$
which is a morphism of correspondence functors because it is induced by the middle vertical morphism
$$L_{E,\CP_E f_R} \longrightarrow L_{E,\CP_E f_R\otimes_{kA} V} \mvirg$$
which is obviously a morphism of correspondence functors.
This defines the morphism $\Phi$ appearing in Proposition~\ref{SERV}.
\fresult

Our goal is to prove that $\Psi:\S_{E,R}\otimes_{kA} V \to S_{E,R,V}$ is an isomorphism.
We prepare the ground with the following lemma, for which we need the full strength of Theorem~\ref{SER-free}, based in turn on Theorem~\ref{dim-SER}.
Since we consider simple modules, we assume that $k$ is a field.

\result{Lemma} \label{Ker-Res-otimes}
Let $k$ be a field, let $(E,R)$ be a finite poset, let $A=\Aut(E,R)$, and let $V$ be a simple $kA$-module.
Let $\alpha \in \S_{E,R}(X)\otimes_{kA}V$ where $X$ is some finite set.
Then $\rho\cdot\alpha=0$ for every $\rho\in k\CC(E,X)$ if and only if $\alpha=0$.
\fresult

\pf
Since $V$ is a simple $kA$-module and $A$ is a finite group, we claim that there exists an injective homomorphism of $kA$-modules $\lambda:V\to kA$.
This follows from the following argument. If $V\dual$ denotes the dual simple module,
there exists a surjective homomorphism $\pi:kA\to V\dual$, which we dualize to obtain an injective homomorphism $\pi\dual:V\to (kA)\dual$.
Now any group algebra is a symmetric algebra, so $(kA)\dual\cong kA$, and this defines the injective homomorphism $\lambda:V\to kA$.\par

If $M$ is a free right $kA$-module, then
$$\Id_M\otimes \lambda: M\otimes_{kA}V \longrightarrow M\otimes_{kA}kA$$
remains injective. This is clear if $M$ is free of rank one and then it follows in general by taking direct sums.
Now we compose with the isomorphism $M\otimes_{kA}kA\cong M$ and we take $M=\S_{E,R}(X)$, which is indeed a free right $kA$-module by Theorem~\ref{SER-free}.
We obtain an injective homomorphism
$$\lambda_X: \S_{E,R}(X)\otimes_{kA}V \longrightarrow \S_{E,R}(X)$$
which is easily seen to define a morphism of correspondence functors
$$\lambda:\S_{E,R}\otimes_{kA}V \longrightarrow \S_{E,R}$$
because we use only the right module structure, whereas correspondences act on the left.\par

For every $\rho\in k\CC(E,X)$, there is a commutative diagram
$$\xymatrix{
\S_{E,R}(X)\otimes_{kA}V \ar[r]^-{\lambda_X} \ar[d]^-{\rho} & \S_{E,R}(X) \ar[d]^-{\rho} \\
\S_{E,R}(E)\otimes_{kA}V \ar[r]^-{\lambda_E} & \S_{E,R}(E)
}$$
Whenever our given element $\alpha \in \S_{E,R}(X)\otimes_{kA}V$ satisfies $\rho\cdot\alpha=0$ for every $\rho\in k\CC(E,X)$,
we also have $\rho\cdot \lambda_X(\alpha)=0$ for every $\rho\in k\CC(E,X)$.
But this implies that $\lambda_X(\alpha)=0$ by Lemma~\ref{Ker-Res}.
Since $\lambda_X$ is injective, $\alpha=0$, as required.
\endpf

Now we come to our main description of simple correspondence functors.

\result{Theorem} \label{iso-SERV}
Let $k$ be a field, let $(E,R)$ be a finite poset, let $A=\Aut(E,R)$, and let $V$ be a simple $kA$-module.
The morphism $\Psi: \S_{E,R}\otimes_{kA}V \to S_{E,R,V}$ is an isomorphism.
\fresult

\pf
For any finite set $X$ and any $\rho\in k\CC(E,X)$, there is a commutative diagram
$$\xymatrix{
\S_{E,R}(X)\otimes_{kA}V \ar[r]^-{\Psi_X} \ar[d]^-{\rho} & S_{E,R,V}(X) \ar[d]^-{\rho} \\
\S_{E,R}(E)\otimes_{kA}V \ar[r]^-{\Psi_E} & S_{E,R,V}(E)
}$$
and $\Psi_E$ is an isomorphism because
$$\S_{E,R}(E)\otimes_{kA}V = \CP_E f_R \otimes_{kA}V =S_{E,R,V}(E) $$
and $\Psi$ is induced by the identity morphism $L_{E,\CP_E f_R}\otimes_{kA} V \to L_{E,\CP_E f_R\otimes_{kA} V}$.\par

Let $\alpha\in \S_{E,R}(X)\otimes_{kA}V$ such that $\Psi_X(\alpha)=0$. Then
$$\Psi_E(\rho\cdot\alpha) = \rho\cdot \Psi_X(\alpha)=0$$
for every $\rho\in k\CC(E,X)$. Since $\Psi_E$ is an isomorphism, we obtain $\rho\cdot\alpha=0$ for every $\rho\in k\CC(E,X)$.
Therefore $\alpha=0$ by Lemma~\ref{Ker-Res-otimes}.
This proves that $\Psi_X$ is injective and we know that it is surjective by construction.
\endpf

We can finally prove one of our main results, namely the determination of the dimension of any evaluation of a simple correspondence functor.
Because of the link with lattices (via the morphism $\Theta_T$), it is convenient to state the result for $R\op$ rather than~$R$.
But this is actually a minor point because $S_{E,R,V}$ is isomorphic to the dual of $S_{E,R\op,V\dual}$ where $V\dual$ is the dual module, by Theorem~9.8 in~\cite{BT3}.

\result{Theorem} \label{dim-SERV}
Let $k$ be a field. Let $(E,R)$ be a poset and let $V$ be a simple left $k\Aut(E,R)$-module.
Let $T$ be a lattice such that $\Irr(T)= (E,R)$ and such that the restriction homomorphism $\Aut(T)\to\Aut(E,R)$ is an isomorphism.
Let $G=E\sqcup\{t\in T\mid t=r^\infty \sigma^\infty(t)\}\subseteq T$ (see Definition~\ref{notation-G} and Lemma~\ref{characterize-G}).\par

For any finite set~$X$, the dimension of $S_{E,R\op,V}(X)$ is given by
$$\dim_kS_{E,R\op,V}(X)=\frac{\dim_kV}{|\Aut(E,R)|}\sum_{i=0}^{|E|}(-1)^i\binom{|E|}{i}(|G|-i)^{|X|} \mpoint$$
\fresult

\pf 
By Theorem~\ref{SER-free}, $\S_{E,R\op}(X)$ is isomorphic to the direct sum of $n_X$ copies of the free right module $k\Aut(E,R)$, for some $n_X\in\N$.
In particular
$$\dim_k\S_{E,R\op}(X)=n_X \, |\Aut(E,R)|\mpoint$$
By Theorem~\ref{iso-SERV}, the simple functor $S_{E,R\op,V}$ is isomorphic to~$\S_{E,R\op}\otimes_{k\Aut(E,R)}V$, using the obvious equality $\Aut(E,R)=\Aut(E,R\op)$.
Thus we obtain
$$S_{E,R\op,V}(X)\cong \S_{E,R\op}(X)\otimes_{k\Aut(E,R)}V\cong n_X \big(k\Aut(E,R)\big)\otimes_{k\Aut(E,R)}V\cong n_XV \mpoint$$
Hence $\dim_k\S_{E,R\op,V}(X)=n_X\dim_kV$.
Therefore
$$\dim_k\S_{E,R\op,V}(X)=\frac{\dim_kV}{|\Aut(E,R)|}\dim_k\S_{E,R\op}(X)\mpoint$$
The result now follows from Theorem~\ref{dim-SER}.  \endpf


\section{Simple modules for the algebra of relations} \label{Section-algebra}

\medskip
\noindent
Let $X$ be a fixed finite set and consider the monoid $\CC(X,X)$ of all relations on~$X$, also known as the monoid of all Boolean matrices of size~$|X|$. As before, write $\CR_X=k\CC(X,X)$ for the algebra of this monoid. Throughout this section, we assume that the base ring~$k$ is a field.
We give the parametrization of all simple modules for the algebra $\CR_X$ and then solve the open problem of giving a formula for their dimension. We also give an explicit description for the action of relations on every simple $\CR_X$-module.\par

We have seen in Theorem~\ref{parametrization} that simple correspondence functors $S_{E,R,V}$ are parametrized by isomorphism classes of triples $(E,R,V)$.
The parametrization of all simple modules for the algebra $\CR_X$ can now be described in terms of simple correspondence functors.

\result{Theorem} \label{simple-RX} Let $X$ be a finite set.
\begin{enumerate}
\item The set of isomorphism classes of simple $\CR_X$-modules is parametrized by the set of isomorphism classes of triples $(E,R,V)$, where $E$ is a finite set with $|E|\leq |X|$, $R$ is an order relation on~$E$, and $V$ is a simple $k\Aut(E,R)$-module.
\item The simple module parametrized by the triple $(E,R,V)$ is $S_{E,R,V}(X)$, where $S_{E,R,V}$ is the simple correspondence functor corresponding to the triple $(E,R,V)$.
\end{enumerate}
\fresult

\pf
We first recall that the evaluation $S(X)$ of a simple correspondence functor~$S$ at a finite set~$X$ is either zero or a simple $\CR_X$-module.
The proof is very easy and appears in Proposition~3.2 of~\cite{We}, or also in Proposition~2.7 of~\cite{BT2}.\par

Conversely, we claim that any simple $\CR_X$-module~$W$ occurs as the evaluation of some simple correspondence functor~$S$, that is, $W\cong S(X)$.
This is Lemma~2.5 of~\cite{BT2} but the proof goes back to the first lemma of~\cite{Bo}.
It also appears in Proposition~3.2 of~\cite{We}, where it is attributed to Green (6.2 in~\cite{Gr}).
This requires to view $\CR_X=k\CC(X,X)$ as a category with a single object~$X$, hence a full subcategory of~$k\CC$.
Proposition~3.2 of~\cite{We} or Proposition~2.7 of~\cite{BT2} also show that the simple correspondence functor~$S$ such that $W\cong S(X)$ is unique up to isomorphism.
All these facts actually hold for the simple representations of any small category.\par

By Theorem~\ref{parametrization}, our simple correspondence functor~$S=S_{E,R,V}$ is parametrized by a triple $(E,R,V)$,
where $E$ is a finite set, $R$ is an order relation on~$E$, and $V$ is a simple $k\Aut(E,R)$-module.
Whenever $W=S_{E,R,V}(X)\neq\{0\}$, we have $|E|\leq|X|$ because $S_{E,R,V}$ vanishes on sets~$Y$ with $|E|>|Y|$ (Proposition~\ref{SERV}).
In order to obtain the parametrization of the statement, we also need to show that $W=S_{E,R,V}(X)$ is nonzero if $|E|\leq |X|$.
This is clear if $|E|=|X|$ because $S_{E,R,V}(E)=T_{R,V}$ is nonzero (see the construction of $S_{E,R,V}$ in Section~\ref{Section-simple-tensor}).
Knowing that $S_{E,R,V}(E)\neq\{0\}$, Corollary~3.7 in~\cite{BT2} asserts precisely that $S_{E,R,V}(X)\neq\{0\}$ if $|E|<|X|$.
This provides the required parametrization and completes the proof.
\endpf

Note that we used in the proof above the non-vanishing of $S_{E,R,V}(X)$ when $|E|<|X|$.
This is a special property of correspondence functors (Corollary~3.7 in~\cite{BT2}) and it may not hold for representations of other small categories.\par

In view of Theorem~\ref{simple-RX}, a formula for the dimension of any simple $\CR_X$-module is now given by Theorem~\ref{dim-SERV}.
More explicitly, we fix a poset $(E,R)$ and a finite lattice $T$ having $(E,R)$ as the full subset of its join-irreducible elements.
We can also choose $T$ such that $\Aut(T)\cong\Aut(E,R)$ by taking for instance $T=\Idown(E,R)$.
We consider the simple $\CR_X$-module $S_{E,R\op,V}(X)$, continuing to use $R\op$ as in Theorem~\ref{dim-SERV}.
We define the subset~$G$ of~$T$ as in Notation~\ref{notation-G} and we write $G=G_{E,R}$ to emphasize its dependence on~$(E,R)$.
Its cardinality $|G|$ only depends on~$(E,R)$, by Corollary~\ref{choice-G}.

\result{Theorem} \label{dim-simple-RX}
With the notation above, the dimension of a simple $\CR_X$-module is given by the formula
$$\dim(S_{E,R\op,V}(X)) = \frac{\dim_kV}{|\Aut(E,R)|} \sum_{i=0}^{|E|}(-1)^i\binom{|E|}{i}(|G_{E,R}|-i)^{|X|}\mpoint$$
\fresult

\pf This is a restatement of Theorem~\ref{dim-SERV}.
\endpf

An explicit description can be given for the action of relations on the simple $\CR_X$-module $S_{E,R\op,V}(X)$.
We define a subset
$$\CB_{E,R,X}=\{\varphi \in T^X \mid E\subseteq \varphi(X) \subseteq G_{E,R} \} \subseteq T^X \mpoint$$
By Theorem~\ref{dim-SER}, the surjective morphism $\Theta_T:F_T\to \S_{E,R\op}$ induces a $k$-module decomposition
$$F_T(X)=k\CB_{E,R,X} \oplus \Ker(\Theta_{T,X}) \mvirg$$
where $k\CB_{E,R,X}$ is the $k$-subspace of~$F_T(X)$ with basis~$\CB_{E,R,X}$.
Thus we have a $k$-module isomorphism
$$\S_{E,R\op}(X) \cong k\CB_{E,R,X} \mpoint$$
The family of subspaces $k\CB_{E,R,X}$ do not form a subfunctor of~$F_T$, but they can be used to describe the evaluations of the functors $\S_{E,R\op}$ and $S_{E,R\op,V}$.\par

We explain a procedure for modifying an element $\varphi\in T^X$ modulo $\Ker(\Theta_{T,X})$ in order to project it in~$k\CB_{E,R,X}$.
In Theorem~\ref{generators-SER}, we introduced an element $u_T\in k(T^T)$ which has the property that, for any $\varphi\in T^X$, the composition $u_T\circ \varphi$ is a $k$-linear combination of maps $f\in T^X$ such that $f(X)\subseteq G_{E,R}$. (Actually, $u_T$ is idempotent, by Theorem~\ref{idempotents-ua}.)
Moreover,
$$u_T\circ \varphi \equiv \varphi \pmod{\Ker(\Theta_{T,X})} \mpoint$$
Let $\pi_{T,X}$ be the $k$-linear idempotent endomorphism of $k(T^X)$ defined by
$$\forall \varphi\in T^X,\;\;\pi_{T,X}(\varphi)=\left\{\begin{array}{ll}\varphi&\hbox{if}\;E \subseteq \varphi(X) \mvirg\\
0&\hbox{otherwise} \mpoint\end{array}\right.$$
By Theorem~\ref{surjection},
$$\pi_{T,X}(\varphi) \equiv \varphi \pmod{\Ker(\Theta_{T,X})} \mpoint$$
Then, for any map $\varphi\in T^X$, we obtain
$$\pi_{T,X}(u_T\circ \varphi) \in k\CB_{E,R,X} \mvirg$$
that is, a $k$-linear combination of maps $f\in T^X$ such that $E\subseteq f(X)\subseteq G_{E,R}$.
Moreover,
$$\pi_{T,X}(u_T\circ \varphi)\equiv \varphi \pmod{\Ker(\Theta_{T,X})} \mpoint$$
Thus if we lift arbitrarily a basis element of~$\S_{E,R\op}(X)$ to $\varphi \in F_T(X)$, we can modify it modulo $\Ker(\Theta_{T,X})$ to obtain an element of~$k\CB_{E,R,X}$. Applying this procedure to the action of a relation $U\in\CC(X,X)$ on an element $\varphi \in k\CB_{E,R,X}$, we obtain
$$U \varphi \equiv \pi_{T,X}(u_T\circ U\varphi) \pmod{\Ker(\Theta_{T,X})} \mvirg$$
and $\pi_{T,X}(u_T\circ U\varphi)$ belongs to~$k\CB_{E,R,X}$.\par

As in Section~\ref{Section-simple-tensor}, we tensor on the right with the $k\Aut(E,R)$-module~$V$, using the right action of $\Aut(E,R)$ on~$\CB_{E,R,X}$ defined by $\varphi\cdot\sigma:= \sigma^{-1}\circ \varphi$ for all $\sigma\in \Aut(E,R)$. By Theorem~\ref{iso-SERV}, we have isomorphisms
$$S_{E,R\op,V}(X) \cong \S_{E,R\op}(X) \otimes_{k\Aut(E,R)} V \cong k\CB_{E,R,X} \otimes_{k\Aut(E,R)} V \mvirg$$
the second isomorphism being only $k$-linear.

This analysis proves the following result, which provides a computational method for describing the action of a relation on a simple $\CR_X$-module.

\result{Theorem} \label{action} Fix the notation above.
\begin{enumerate}
\item $S_{E,R\op,V}(X)\cong k\CB_{E,R,X} \otimes_{k\Aut(E,R)} V$ as $k$-vector spaces.
\item Transporting the action of relations via this isomorphism, the action of a relation $U\in\CC(X,X)$ on an element
$$\varphi\otimes v \in k\CB_{E,R,X} \otimes_{k\Aut(E,R)} V \,, \qquad(\varphi\in\CB_{E,R,X}, \; v\in V)$$
is given by
$$U\cdot (\varphi\otimes v)=\pi_{T,X}(u_T\circ U\varphi) \otimes v \mpoint$$
\end{enumerate}
\fresult

Our last result gives the dimension of the Jacobson radical $J(\CR_X)$ of the $k$-algebra~$\CR_X$.
We assume for simplicity that the field $k$ has characteristic zero. 

\result{Theorem} \label{Jacobson}
Assume that $k$ is a field of characteristic zero.
Let $J(\CR_X)$ be the Jacobson radical of the $k$-algebra~$\CR_X$ and let $n=|X|$. Then
$$\dim J(\CR_X) = 2^{n^2} - \sum_{e=0}^n \sum_R \displaystyle\frac{1}{|\Aut(E,R)|}
\Big(\sum_{i=0}^e (-1)^i \binom{e}{i}(|G_{E,R}|-i)^{n}  \Big)^2 \mvirg$$
where $R$ runs over a set of representatives of $\Sigma_e$-conjugacy classes of order relations on the set~$E=\{1,\ldots,e\}$ and $E=\emptyset$ if $e=0$.
The integer $|G_{E,R}|$ is the cardinality of the set $G_{E,R}$ defined in Notation~\ref{notation-G}.
\fresult

\pf
Since $k$ has characteristic zero, the semi-simple algebra $\CR_X/J(\CR_X)$ is separable, that is, it remains semi-simple after scalar extension to an algebraic closure $\overline k$ of~$k$. In other words, $\dim J(\CR_X)$ does not change after this scalar extension. Therefore, we can assume that $k=\overline k$.\par

By Theorem~\ref{simple-RX}, every simple $\CR_X$-module has the form $S_{E,R,V}(X)$ with $|E|\leq|X|$, where $S_{E,R,V}$ is the simple correspondence functor parametrized by the triple $(E,R,V)$. In order to have a parametrization, we take $E=\{1,\ldots,e\}$ with ${0\leq e\leq n}$, we take $R$ in a set of representatives as in the statement, and finally we take $V$ in a set of representatives of isomorphism classes of simple $k\Aut(E,R)$-modules.\par

Since the endomorphism algebra of a simple module is isomorphic to~$k$, by Schur's lemma and the assumption that $k$ is algebraically closed, the dimension of the semi-simple algebra $\CR_X/J(\CR_X)$ is equal to the sum of the squares of the dimensions of all simple modules, by Wedderburn's theorem. From the formula for the dimension of simple modules, we obtain
$$\begin{array}{rcl}
\dim \big(\CR_X/J(\CR_X)\big) &=& \displaystyle \sum_{E,R,V} \big(\dim S_{E,R,V}(X)\big)^2 = \displaystyle \sum_{E,R,V} \big(\dim S_{E,R\op,V}(X)\big)^2\\
&=& \displaystyle \sum_{E,R,V} \Big(\frac{\dim V}{|\Aut(E,R)|}\Big)^2 \Big(\sum_{i=0}^{|E|}(-1)^i\binom{|E|}{i}(|G_{E,R}|-i)^{|X|}\Big)^2 \\
&=& \displaystyle \sum_{e=0}^n \sum_R \Big(\sum_V \frac{(\dim V)^2}{|\Aut(E,R)|^2}\Big)
\Big(\sum_{i=0}^e (-1)^i \binom{e}{i}(|G_{E,R}|-i)^{n}\Big)^2 \\
&=& \displaystyle \sum_{e=0}^n \sum_R \frac{1}{|\Aut(E,R)|}
\Big(\sum_{i=0}^e (-1)^i \binom{e}{i}(|G_{E,R}|-i)^{n}\Big)^2
\mvirg\end{array}$$
because $\displaystyle\sum_V(\dim V)^2=\dim (k\Aut(E,R)) =|\Aut(E,R)|$, by semi-simplicity of the group algebra in characteristic zero (Maschke's theorem). Now
$$\dim J(\CR_X) = \dim \CR_X- \dim \big(\CR_X/J(\CR_X)\big) = 2^{n^2} - \dim \big(\CR_X/J(\CR_X)\big)$$
and the result follows.
\endpf

If $k$ is an algebraically closed field of prime characteristic~$p$, the formula has to be modified in a straightforward manner, in order to take into account the Jacobson radical of $k\Aut(E,R)$.
Then it seems likely that the same formula holds over any field of characteristic~$p$ (that is, $\CR_X/J(\CR_X)$ is likely to be a separable algebra), but we leave this question open.


\section{Examples} \label{Section-examples}

\medskip
\noindent
We state here without proofs a list of examples. For simplicity, we assume that the base ring~$k$ is a field (but many results actually remain true over an arbitrary commutative ring~$k$).
We first describe a few small cases for modules over the algebra $\CR_X$, using the notation of Section~\ref{Section-algebra}.
Then we give the decomposition of the functors~$F_T$ associated to some particular lattices~$T$.

\subsect{Example} \label{empty}
Let $X=\emptyset$. There is a single relation on~$\emptyset$, namely~$\emptyset$, and $\CR_\emptyset\cong k$. Then
$S_{\emptyset,\emptyset,k}(\emptyset)\cong k\CB_{\emptyset,\emptyset,\emptyset} \otimes_k k\cong k$
and the unique relation $\emptyset$ acts as the identity on~$k$.

\subsect{Example} \label{one}
Let $X=\{1\}$. There are 2 relations on~$\{1\}$, namely $\emptyset$ and $\Delta_{\{1\}}$.\par

For $E=\emptyset$, we get $S_{\emptyset,\emptyset,k}(\{1\})\cong k\CB_{\emptyset,\emptyset,\{1\}} \otimes_k k\cong k$
and both relations act as the identity on~$k$.\par

For $E=\{1\}$, we obtain $S_{\{1\},\Delta_{\{1\}},k}(\{1\})\cong k\CB_{\{1\},\Delta_{\{1\}},\{1\}} \otimes_k k \cong k$,
the relation $\emptyset$ acts by zero, while $\Delta_{\{1\}}$ acts as the identity.\par

Moroever, $\CR_{\{1\}}\cong k\times k$ is a semi-simple algebra.

\subsect{Example} \label{two}
Let $X=\{1,2\}$. There are $2^4=16$ relations on~$\{1,2\}$, so $\CR_{\{1,2\}}$ has dimension~16.\par

For $E=\emptyset$, we get a simple $\CR_{\{1,2\}}$-module $S_{\emptyset,\emptyset,k}(\{1,2\})$ of dimension~1.\par

For $E=\{1\}$, we get a simple $\CR_{\{1,2\}}$-module $S_{\{1\},\Delta_{\{1\}},k}(\{1,2\})$ of dimension~3.\par

For $E=\{1,2\}$, there are two essential relations up to conjugacy, namely the equality relation~$\Delta_{\{1,2\}}$ and the usual total order $\rm{tot}$.
Moreover, $\Aut(\{1,2\},\Delta_{\{1,2\}})$ is a group of order~2, with two simple modules $k_+$ and $k_-$ (assuming that the characteristic of $k$ is not~2).
Therefore, we obtain two simple $\CR_{\{1,2\}}$-modules of dimension~1
$$\begin{array}{rcl}
S_{\{1,2\},\Delta_{\{1,2\}},k_+}(\{1,2\}) &\cong& k\CB_{\{1,2\},\Delta_{\{1,2\}},\{1,2\}} \otimes_{kC_2} k_+ \mvirg \\
S_{\{1,2\},\Delta_{\{1,2\}},k_-}(\{1,2\}) &\cong& k\CB_{\{1,2\},\Delta_{\{1,2\}},\{1,2\}} \otimes_{kC_2} k_- \mvirg
\end{array}
$$
For the other relation $\rm{tot}$, we obtain a simple $\CR_{\{1,2\}}$-module of dimension~2
$$S_{\{1,2\},{\rm tot},k}(\{1,2\}) \cong k\CB_{\{1,2\},\rm{tot},\{1,2\}} \otimes_k k\cong k\CB_{\{1,2\},\rm{tot},\{1,2\}}
\mvirg$$
Altogether, there are 5 simple $\CR_{\{1,2\}}$-modules and the Jacobson radical has dimension~0.
Therefore $\CR_{\{1,2\}}$ is semi-simple (provided the characteristic of $k$ is not~2).

\def\posetVop{^{\mathop{\mathop{\displaystyle\bullet\;\bullet}^{\displaystyle/\backslash}\limits}^{\raisebox{-1ex}{$\displaystyle\bullet$}}\limits}}
\def\posetV{\raisebox{-2ex}{$^{\,\mathop{\mathop{\displaystyle\bullet\;\bullet}_{\displaystyle\backslash/}\limits}_{^{\displaystyle\bullet}}\limits}$}}
\def\ipoint{\mathop{\rule{.1ex}{.9ex}\hspace{.2ex}}^{_{\displaystyle\bullet}}_{^{\displaystyle\bullet}}\limits{\displaystyle\bullet}}
\def\totdeux{\mathop{\rule{.1ex}{.9ex}\hspace{.2ex}}\limits^{_{\displaystyle\bullet}}_{^{\displaystyle\bullet}}}
\def\tottrois{\mathop{\totdeux}\limits^{\mathop{\rule{.1ex}{.9ex}}\limits^{\displaystyle\bullet}}}
\subsect{Example} \label{three}
For $|X|= 3$, the algebra $\CR_X$ is not  semi-simple. The dimension of the Jacobson radical of~$\CR_X$ is equal to~42, using either the computer software~\cite{GAP4} or the computer calculations obtained in~\cite{Br}.
According to Theorem~\ref{Jacobson}, this value can be recovered directly as follows~:
\medskip

$$\begin{array}{|c|c|c|c|c|r|}
\hline
{\rm Size}\;e&{\rm Poset}(E,R)&|\Aut(E,R)|&|G_{E,R}|&\sum\limits_{i=0}^e(-1)^i\binom{e}{i}(|G_{E,R}|-i)^3&{\rm total}\\
\hline
0&\emptyset&1&1&1&1\\
\hline
1&\bullet&1&2&7&49\\
\hline
2&\bullet\bullet&2&4&18&162\\
&\rule{0ex}{4ex}\totdeux&1&3&12&144\\
\hline
3&\bullet\bullet\bullet&6&5&6&6\\
&\rule{0ex}{4ex}\ipoint&1&5&6&36\\
&\rule{0ex}{5ex}\posetV&2&5&6&18\\
&\rule{0ex}{5ex}\raisebox{-2ex}{$\posetVop$}&2&5&6&18\\
&\tottrois&1&6&6&36\\
\hline
\end{array}
$$
In this case, the algebra $\CR_X$ has dimension $2^{3^2}=512$. The sum of the last column of this table is equal to 470, so we recover the dimension of the radical $42=512-470$.

\medskip

\subsect{Example} \label{more}
For $|X|= 4$, the algebra $\CR_X$ has dimension $2^{4^2}=65,\!536$. The direct computation of the radical of such a big algebra seems out of reach of usual computers. However, using the formula of Theorem~\ref{Jacobson} and the structure of the 16 posets of cardinality 4, one can show by hand that the radical of $\CR_X$ has dimension $32,\!616$.\par

For larger values of $n=|X|$, a computer calculation using Theorem~\ref{Jacobson} yields the following values for the dimension of $J(\CR_X)$~:
$$\begin{array}{|c|c|c|c|}
\hline
n=5&n=6&n=7&n=8\\
\hline
 29,\!446,\!050 & 67,\!860,\!904,\!320 & 562,\!649,\!705,\!679,\!642 & 18,\!446,\!568,\!932,\!288,\!588,\!616 \\
\hline
\end{array}
$$

\bigskip
We now move to examples of fundamental functors and functors associated to lattices.

\subsect{Example} \label{forest}
There are many examples of fundamental functors $\S_{E,R}$ for which the set $G$ is the whole of~$T$, for instance when $T=\Lambda E$.
In all such cases, we have $F_T/H_T\cong\S_{E,R}$.
Moreover, in many such cases, $\Aut(E,R)$ is the trivial group. Take for instance $(E,R)$ to be a disjoint union of trees with branches of different length.
In any such case, $F_T/H_T\cong\S_{E,R}\cong S_{E,R,k}$ is simple, provided $k$ is a field.

\bigskip
Our next purpose is to decompose the functor $F_T$ for some small lattices~$T$.
In order to use an inductive process, we use inclusions $A\to T$ where $A$ is a distributive sublattice.
We could as well use surjective morphisms $T\to A$, as in Section~10 of~\cite{BT3}, but the following general result shows that it does not matter.

\result{Lemma} \label{split}
Let $T$ and $A$ be finite lattices and assume that $A$ is distributive.
Let $\sigma:A\to T$ be an injective join-preserving map.
Then there is a surjective join-preserving map $\pi:T\to A$ such that $\pi\sigma=\Id_A$.
\fresult

\pf
We define $\pi(t)=\displaystyle\bigwedge_{\substack{a\in A\\ \sigma(a)\geq t}}a$.
Then $\pi$ is order-preserving and therefore $\pi(t_1)\vee \pi(t_2)\leq \pi(t_1\vee t_2)$ for any $t_1,t_2\in T$.
Now we have
$$\pi(t_1)\vee \pi(t_2)
=\big( \bigwedge_{\substack{a_1\in A\\ \sigma(a_1)\geq t_1}}a_1 \big) \vee \big( \bigwedge_{\substack{a_2\in A\\ \sigma(a_2)\geq t_2}}a_2 \big)
=\bigwedge_{\substack{a_1,a_2\in A \\ \sigma(a_1)\geq t_1 \\ \sigma(a_2)\geq t_2}} (a_1\vee a_2)$$
by distributivity of~$A$.
For any such pair $(a_1,a_2)$, the join $a_1\vee a_2$ belongs to the set $\{a\in A \mid \sigma(a)\geq t_1\vee t_2 \}$ and therefore
$$\bigwedge_{\substack{a_1,a_2\in A \\ \sigma(a_1)\geq t_1 \\ \sigma(a_2)\geq t_2}} (a_1\vee a_2)
\;\geq \bigwedge_{\substack{a\in A\\ \sigma(a)\geq t_1\vee t_2}}a \; = \; \pi(t_1\vee t_2) \mpoint$$
The equality $\pi(t_1)\vee \pi(t_2)= \pi(t_1\vee t_2)$ follows.\par

If $\sigma(a_1)\leq \sigma(a_2)$, then $\sigma(a_1\vee a_2)=\sigma(a_1)\vee\sigma(a_2)=\sigma(a_2)$,
and therefore $a_1\vee a_2=a_2$ by injectivity of~$\sigma$, i.e. $a_1\leq a_2$. It follows from this observation that, for any $b\in A$,
$$\pi\sigma(b) =\bigwedge_{\substack{a\in A\\ \sigma(a)\geq \sigma(b)}}a \;=\bigwedge_{\substack{a\in A\\ a\geq b}}a = \; b \mvirg$$
hence $\pi\sigma=\Id_A$.
\endpf

The property of Lemma~\ref{split} is reflected in the fact that the morphism $F_A\to F_T$ induced by~$\sigma$ must split,
because the functor $F_A$ is projective by Theorem~4.12 in~\cite{BT3} and injective by Theorem~10.6 in~\cite{BT2}.

\subsect{Example} \label{lozenge}
Let $T=\lozenge$ be the lozenge, in other words the lattice of subsets of a set of cardinality 2~:
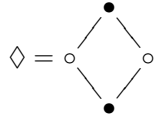
$$\lozenge=\vcenter{\xymatrix@R=2ex@C=2ex{&\pointplein\ar@{-}[ld]\ar@{-}[rd]&\\
\pointcreux\ar@{-}[rd]&&\pointcreux\ar@{-}[ld]\\
&\pointplein&
}
}
$$
By Theorem~11.12 in~\cite{BT3}, for any finite lattice~$T$, we can split off from~$F_T$ simple functors $\S_n:=\S_{\sou{n},\tot} \cong S_{\sou{n},\tot,k}$ corresponding to all totally ordered sequences $\widehat 0\leq d_0<d_1<\ldots<d_n=\widehat 1$ in~$T$.
In the case of $F_\lozenge$, we obtain
$$F_\lozenge\cong \S_0\oplus 3\S_1\oplus 2\S_2\oplus L$$
for some subfunctor $L$.
\def\bcirc{{\scriptscriptstyle\circ\circ}}
We know that $F_\lozenge$ maps surjectively onto the fundamental functor $\S_{\bcirc}$ associated to the (opposite) poset of irreducible elements of~$\lozenge$, that is, a set of cardinality~2 ordered by the equality relation. Moreover, all the factors $\S_n$ lie in the subfunctor $H_\lozenge$, because no totally ordered subset contains the two irreducible elements of $\lozenge$ (figured with an empty circle in the above picture).
Therefore $L$ maps surjectively onto~$\S_{\bcirc}$.

We can evaluate $F_\lozenge$ at a set $X$ of cardinality~$x$, and take dimensions over~$k$. By Theorem~\ref{dim-SER}, we obtain
$$4^x=1^x+3(2^x-1^x)+2(3^x-2\cdot 2^x+1^x)+\dim_kL(X)\mpoint$$ 
It follows that 
\def\bcirc{{\scriptscriptstyle\circ\circ}}
$$\dim_kL(X)=4^x-2\cdot 3^x+2^x\mpoint$$
Now we apply Theorem~\ref{dim-SER} to the fundamental functor $\S_{\bcirc}$.
The set $G$ is the whole of~$T$ in this case, so $\dim_k\S_{\bcirc}(X)=4^x-2\cdot 3^x+2^x$. \par

Since $L$ maps surjectively onto $\S_{\bcirc}$ and $\dim_kL(X)=\dim_k\S_{\bcirc}(X)$ for any finite set $X$, this surjection is an isomorphism. Hence
$$F_\lozenge\cong \S_0\oplus 3\S_1\oplus 2\S_2\oplus \S_{\bcirc}\mpoint$$
Since the lattice $\lozenge$ is distributive, $F_\lozenge$ is projective by Theorem~4.12 in~\cite{BT3}.
It follows that each summand is a projective object in the category~$\CF_k$ of correspondence functors.

\subsect{Example} \label{equality-3}
Let $T$ be the following lattice~:
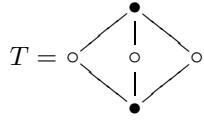
$$T=\vcenter{\xymatrix@R=2ex@C=4ex{&\pointplein\ar@{-}[ld]\ar@{-}[rd]\ar@{-}[d]&\\
\pointcreux\ar@{-}[rd]&\pointcreux\ar@{-}[d]&\pointcreux\ar@{-}[ld]\\
&\pointplein&
}
}
$$
As in the previous example, $F_T$ admits a direct summand isomorphic to 
$$\S_0\oplus 4\S_1\oplus 3\S_2\mpoint$$
Moreover, there are three obvious sublattices of~$T$ isomorphic to~$\lozenge$, which provide three direct summands of~$F_T$ isomorphic to~$\S_{\bcirc}$. Thus we have a decomposition
$$F_T\cong \S_0\oplus 4\S_1\oplus 3\S_2\oplus3\S_{\bcirc}\oplus M$$
for some subfunctor $M$ of~$F_T$. 
\def\ccirc{{\scriptscriptstyle\circ\circ\circ}}
Using arguments similar to those of the previous example, we get
$$F_T\cong \S_0\oplus 4\S_1\oplus 3\S_2\oplus3\S_{\bcirc}\oplus \S_{\ccirc}\mpoint$$
All the summands in this decomposition of~$F_T$, except possibly $\S_{\ccirc}$, are projective functors. Since the lattice $T$ is not distributive, the functor $F_T$ is not projective (Theorem~4.12 in~\cite{BT3}), thus $\S_{\ccirc}$ is actually not projective either.

\subsect{Example} \label{diamond}
Let $D$ be the following lattice~:
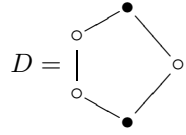
$$D=\vcenter{\xymatrix@R=.1ex@C=3ex{&\pointplein\ar@{-}[ld]\ar@{-}[rdd]&\\
\pointcreux\ar@{-}[dd]&&\\
&&\pointcreux\ar@{-}[ldd]\\
\pointcreux\ar@{-}[rd]&&\\
&\pointplein&
}
}
$$
As before, we know that $F_T$ admits a direct summand isomorphic to a direct sum $\S_0\oplus 4\S_1\oplus 4\S_2\oplus\S_3$.
Moreover, there are two inclusions
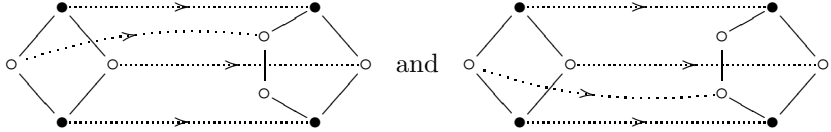
\def\marg[#1]{\ar@{-}[#1]|-{\object@{<}}}
\def\mard[#1]{\ar@{.}[#1]|-{\object@{>}}}
\def\mardb[#1]{\ar@{.}@/^/[#1]|-{\object@{>}}}
\def\mardd[#1]{\ar@{.}@/_/[#1]|-{\object@{>}}}
\def\marb[#1]{\ar@{-}[#1]|{\object+{  }}}
\def\ipoint{\mathop{\rule{.1ex}{.9ex}\hspace{.2ex}}^{_\circ}_{^\circ}\limits{\scriptscriptstyle\circ}}
\def\iipoint{\;\mathop{\rule{.1ex}{.7ex}\hspace{.1ex}}^{_\circ}_{^\circ}\limits{\scriptscriptstyle\circ}}
$$\vcenter{\xymatrix@R=.1ex@C=3ex{
&\pointplein\ar@{-}[ldd]\ar@{-}[rdd]\mard[rrrrr]&&&&&\pointplein\ar@{-}[ld]\ar@{-}[rdd]&\\
&&&&&\pointcreux\ar@{-}[dd]&&\\
\pointcreux\ar@{-}[rdd]\mardb[urrrrr]&&\pointcreux\ar@{-}[ldd]\mard[rrrrr]&&&&&\pointcreux\ar@{-}[ldd]\\
&&&&&\pointcreux\ar@{-}[rd]&&\\
&\pointplein\mard[rrrrr]&&&&&\pointplein&
}
}
\;\;\;\hbox{and}\;\;\;
\vcenter{\xymatrix@R=.1ex@C=3ex{
&\pointplein\ar@{-}[ldd]\ar@{-}[rdd]\mard[rrrrr]&&&&&\pointplein\ar@{-}[ld]\ar@{-}[rdd]&\\
&&&&&\pointcreux\ar@{-}[dd]&&\\
\pointcreux\ar@{-}[rdd]\mardd[drrrrr]&&\pointcreux\ar@{-}[ldd]\mard[rrrrr]&&&&&\pointcreux\ar@{-}[ldd]\\
&&&&&\pointcreux\ar@{-}[rd]&&\\
&\pointplein\mard[rrrrr]&&&&&\pointplein&
}
}
$$
of the lattice $\lozenge$ into $D$, which yield two direct summands of $F_D$ isomorphic to~$\S_{\bcirc}$.
So there is a decomposition
$$F_D\cong \S_0\oplus 4\S_1\oplus 4\S_2\oplus\S_3\oplus 2\S_{\bcirc}\oplus N$$
for a suitable subfunctor $N$ of~$F_D$. As in the previous examples, the subfunctor~$N$ maps surjectively onto the fundamental functor $\S_{\iipoint}$ associated to the (opposite) poset $\ipoint$ of irreducible elements of~$D$. Computing dimensions, we obtain $N\cong \S_{\iipoint}$, and therefore
$$F_D\cong \S_0\oplus 4\S_1\oplus 4\S_2\oplus\S_3\oplus 2\S_{\bcirc}\oplus\S_{\iipoint}\mpoint$$
Again $D$ is not distributive, so that $F_D$ is not projective. Thus the functor $\S_{\iipoint}$ is not projective either.\par

Actually, the lattice $D$ and the lattice $T$ of the previous example are the smallest non-distributive lattices
and they are used for the well-known characterization of distributive lattices (see Theorem~4.7 in~\cite{Ro}).

\subsect{Example} \label{lattice-V}
Let $C$ be the following lattice~:
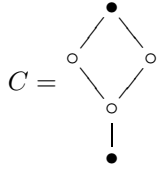
$$C=\vcenter{\xymatrix@R=2ex@C=2ex{
&\pointplein\ar@{-}[ld]\ar@{-}[rd]&\\
\pointcreux\ar@{-}[rd]&&\pointcreux\ar@{-}[ld]\\
&\pointcreux\ar@{-}[d]&\\
&\pointplein&
}
}
$$
Again, we know that $F_C$ admits a direct summand isomorphic to a direct sum $\S_0\oplus 4\S_1\oplus 5\S_2\oplus2\S_3$.
Moreover, the inclusion
$$
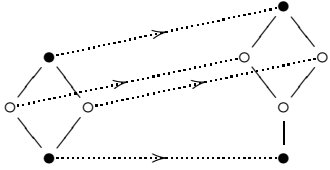
\xymatrix@R=2ex@C=2ex{
&&&&&&&\pointplein\ar@{-}[ld]\ar@{-}[rd]&\\
&\pointplein\ar@{-}[ld]\ar@{-}[rd]\mard[urrrrrr]&&&&&\pointcreux\ar@{-}[rd]&&\pointcreux\ar@{-}[ld]\\
\pointcreux\ar@{-}[rd]\mard[urrrrrr]&&\pointcreux\ar@{-}[ld]\mard[urrrrrr]&&&&&\pointcreux\ar@{-}[d]&\\
&\pointplein\mard[rrrrrr]&&&&&&\pointplein&
}
$$
of $\lozenge$ in $C$ yields a direct summand of $F_C$ isomorphic to~$\S_{\bcirc}$. So there is a decomposition
$$F_C\cong \S_0\oplus 4\S_1\oplus 5\S_2\oplus2\S_3\oplus\S_{\bcirc}\oplus Q$$
for some direct summand $Q$ of~$F_C$.
\def\posetVop{^{\mathop{\mathop{\scriptscriptstyle\circ\;\circ}^{\scriptscriptstyle/\backslash}\limits}^{\raisebox{-1ex}{$\scriptscriptstyle\circ$}}\limits}}

Now $F_C$ maps surjectively onto the fundamental functor $\S_{\posetVop}$ associated to the opposite poset of its irreducible elements, and arguments as before yield an isomorphism $Q\cong \S_{\posetVop}$, hence a decomposition
$$F_C\cong \S_0\oplus 4\S_1\oplus 5\S_2\oplus2\S_3\oplus\S_{\bcirc}\oplus\S_{\posetVop}\mpoint$$
Since $C$ is distributive, $F_C$ is projective and we conclude that $\S_{\posetVop}$ is projective.
Taking dual functors corresponds to taking opposite lattices (see Theorem~8.9 and Remark~9.7 in~\cite{BT3}), so we get a decomposition
\def\posetV{^{\,\mathop{\mathop{\scriptscriptstyle\circ\;\circ}_{\scriptscriptstyle\backslash/}\limits}_{^\circ}\limits}}
$$F_{C\op}\cong\S_0\oplus 4\S_1\oplus 5\S_2\oplus2\S_3\oplus\S_{\bcirc}\oplus\S_{\posetV}\mpoint$$
Therefore $\S_{\posetV}$ is also projective.

\subsect{Example} \label{2-times-3}
Let $P$ be the following lattice~:
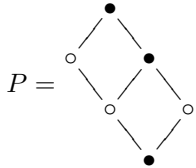
$$P=\vcenter{\xymatrix@R=2ex@C=2ex{&\pointplein\ar@{-}[ld]\ar@{-}[rd]&&\\
\pointcreux\ar@{-}[rd]&&\pointplein\ar@{-}[ld]\ar@{-}[rd]&\\
&\pointcreux\ar@{-}[rd]&&\pointcreux\ar@{-}[ld]\\
&&\pointplein&
}
}
$$
that is, the direct product of a totally ordered lattice of cardinality~3 with a totally ordered lattice of cardinality~2.\par
We know that $F_P$ admits a direct summand isomorphic to $\S_0\oplus 5\S_1\oplus 7\S_2\oplus3\S_3$ and the inclusions
$$
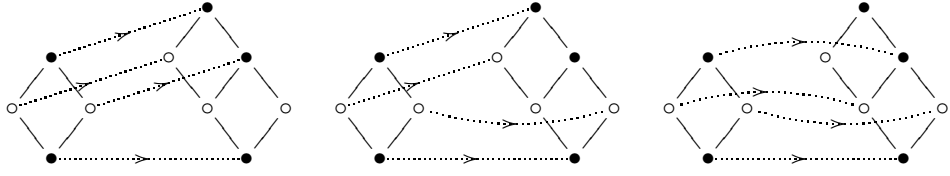
\vcenter{\xymatrix@R=2ex@C=2ex{
&&&&&&\pointplein\ar@{-}[ld]\ar@{-}[rd]&&\\
&&\pointplein\ar@{-}[ld]\ar@{-}[rd]\mard[urrrr]&&&\pointcreux\ar@{-}[rd]&&\pointplein\ar@{-}[ld]\ar@{-}[rd]&\\
&\pointcreux\ar@{-}[rd]\mard[urrrr]&&\pointcreux\ar@{-}[ld]\mard[urrrr]&&&\pointcreux\ar@{-}[rd]&&\pointcreux\ar@{-}[ld]\\
&&\pointplein\mard[rrrrr]&&&&&\pointplein&
}
}
\vcenter{\xymatrix@R=2ex@C=2ex{
&&&&&&\pointplein\ar@{-}[ld]\ar@{-}[rd]&&\\
&&\pointplein\ar@{-}[ld]\ar@{-}[rd]\mard[urrrr]&&&\pointcreux\ar@{-}[rd]&&\pointplein\ar@{-}[ld]\ar@{-}[rd]&\\
&\pointcreux\ar@{-}[rd]\mard[urrrr]&&\pointcreux\ar@{-}[ld]\mardd[rrrrr]&&&\pointcreux\ar@{-}[rd]&&\pointcreux\ar@{-}[ld]\\
&&\pointplein\mard[rrrrr]&&&&&\pointplein&
}
}
\vcenter{\xymatrix@R=2ex@C=2ex{
&&&&&&\pointplein\ar@{-}[ld]\ar@{-}[rd]&&\\
&&\pointplein\ar@{-}[ld]\ar@{-}[rd]\mardb[rrrrr]&&&\pointcreux\ar@{-}[rd]&&\pointplein\ar@{-}[ld]\ar@{-}[rd]&\\
&\pointcreux\ar@{-}[rd]\mardb[rrrrr]&&\pointcreux\ar@{-}[ld]\mardd[rrrrr]&&&\pointcreux\ar@{-}[rd]&&\pointcreux\ar@{-}[ld]\\
&&\pointplein\mard[rrrrr]&&&&&\pointplein&
}
}
$$
of $\lozenge$ in $P$ yield 3 direct summands of $F_P$ isomorphic to $\S_{\bcirc}$. 
Moreover, the inclusions
$$\vcenter
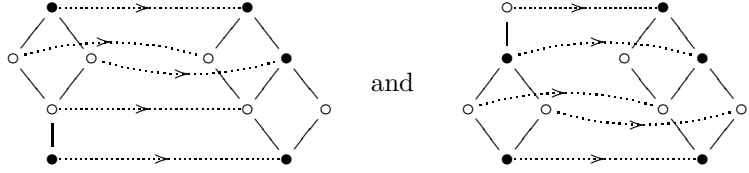
{\xymatrix@R=2ex@C=2ex{
&\pointplein\ar@{-}[ld]\ar@{-}[rd]\mard[rrrrr]&&&&&\pointplein\ar@{-}[ld]\ar@{-}[rd]&&\\
\pointcreux\ar@{-}[rd]\mardb[rrrrr]&&\pointcreux\ar@{-}[ld]\mardd[rrrrr]&&&\pointcreux\ar@{-}[rd]&&\pointplein\ar@{-}[ld]\ar@{-}[rd]&\\
&\pointcreux\ar@{-}[d]\mard[rrrrr]&&&&&\pointcreux\ar@{-}[rd]&&\pointcreux\ar@{-}[ld]\\
&\pointplein\mard[rrrrrr]&&&&&&\pointplein&
}
}\;\;\;\;\;\hbox{and}\;
\vcenter{\xymatrix@R=2ex@C=2ex{
&&\pointcreux\ar@{-}[d]\mard[rrrr]&&&&\pointplein\ar@{-}[ld]\ar@{-}[rd]&&\\
&&\pointplein\ar@{-}[ld]\ar@{-}[rd]\mardb[rrrrr]&&&\pointcreux\ar@{-}[rd]&&\pointplein\ar@{-}[ld]\ar@{-}[rd]&\\
&\pointcreux\ar@{-}[rd]\mardb[rrrrr]&&\pointcreux\ar@{-}[ld]\mardd[rrrrr]&&&\pointcreux\ar@{-}[rd]&&\pointcreux\ar@{-}[ld]\\
&&\pointplein\mard[rrrrr]&&&&&\pointplein&
}
}
$$
of $C$ and $C\op$ in $P$ yield direct summands $\S_{\posetV}$ and~$\S_{\posetVop}$ of~$F_P$, hence there is a direct summand $U$ of~$F_P$ such that
\def\ivague{\mathop{/\backslash/}^{\;\,\circ\;\;\circ}_{\circ\;\;\circ\;\,}\limits}
$$F_P\cong \S_0\oplus 5\S_1\oplus 7\S_2\oplus3\S_3\oplus 3\S_{\bcirc}\oplus \S_{\posetV}\oplus \S_{\posetVop}\oplus U\mpoint$$
Since the lattice $P$ is distributive, the functor $F_P$ is projective, hence $U$ is projective. Now $F_P$ maps surjectively onto the fundamental functor~$\S_{\iipoint}$, and $H_P$ is contained in the kernel of this surjection. It follows that $U$ maps surjectively onto $\S_{\iipoint}$, which is a simple functor, as $k$ is a field and the poset $\ipoint$ has no nontrivial automorphisms.\par

A more involved analysis shows that $U$ is indecomposable and is a projective cover of the simple functor~$\S_{\iipoint}$.
Moreover, one can show that the functor $U$ is uniserial, with a filtration 
$$\xymatrix{0\ar@{-}@/_3ex/[rr]_-{\S_{\iipoint}}&\subset&W\ar@{-}@/_3ex/[rr]_-{\S_{\ivague}}&\subset&V\ar@{-}@/_3ex/[rr]_-{\S_{\iipoint}}&\subset&U\mvirg}
$$
where $W\cong U/V\cong \S_{\iipoint}$, and $V/W$ is isomorphic to the simple functor $\S_{^{\ivague}}$ associated to the poset 
$\vcenter{\xymatrix@R=1ex@C=1ex{
&\pointcreux\ar@{-}[ld]\ar@{-}[rd]&&\pointcreux\ar@{-}[ld]\\
\pointcreux&&\pointcreux&
}
}$
of cardinality~4. An easy consequence of this is that
$${\rm Ext}^1_{\CF_k}(\S_{\iipoint},\S_{\ivague})\cong {\rm Ext}^1_{\CF_k}(\S_{\ivague},\S_{\iipoint})\cong k\mpoint$$


\bigskip
\noindent
Serge Bouc, CNRS-LAMFA, Universit\'e de Picardie - Jules Verne,\\
33, rue St Leu, F-80039 Amiens Cedex~1, France.\\
{\tt serge.bouc@u-picardie.fr}

\medskip
\noindent
Jacques Th\'evenaz, Section de math\'ematiques, EPFL, \\
Station~8, CH-1015 Lausanne, Switzerland.\\
{\tt Jacques.Thevenaz@epfl.ch}

\end{document}